\documentstyle[12pt]{article}

\def\IR{\hbox{\rm R \kern-1.1em \vrule depth 0ex height 1.55ex width .05em
\kern.41em}}

\def\complexes{\hbox{\rm C \kern-.80em \vrule depth 0ex height 1.5ex width
.05em \kern.41em}}
\def\IC{\complexes}

\def\core{{\rm core}}
\def\ord{{\rm ord}}
\def\sp{{\rm sp\!}}
\def\span{{\rm span}\,}
\def\supp{{\rm supp}\,}
\def\rbd{{\rm rbd}\,}
\def\relint{{\rm relint}\,}
\def\rank{{\rm rank}\,}
\def\int{{\rm int}\,}
\def\epi{{\rm epi}\,}
\def\conv{{\rm conv}\,}
\def\order{{\rm order}\,}
\def\pos{{\rm pos}\,}

\def\lline{
\unitlength 1mm
\begin{picture}(10,4)
\put(1,0){\line(1,0){8}}
\end{picture}
}

\begin{document}
\setcounter{page}{0}

{\baselineskip 0.6 cm
\begin{center}
{\Large ON THE INVARIANT FACES ASSOCIATED}
\end{center}
\begin{center}
{\Large WITH A CONE-PRESERVING MAP}
\end{center}
\vskip 40 pt
\begin{center}
{\large Bit-Shun Tam}
\end{center}
\par\parindent=120 pt
Department of Mathematics

Tamkang University

Tamsui,
Taiwan 25137,
ROC

(e-mail: bsm01{$@$}mail.tku.edu.tw;

fax: 2-6209916)

\vskip 12 pt
\begin{center}
and
\end{center}
\vskip 12 pt
\begin{center}
{\large Hans Schneider}
\end{center}
\par\parindent=120 pt
Department of Mathematics

University of Wisconsin

Madison,
Wisconsin 53706,
USA

(e-mail: hans{$@$}math.wisc.edu;

fax: 608-263-8891)

\vskip 30 pt
\begin{center}
22 October 1997
\end{center}

\par\parindent=0 pt
\vskip 78 pt

Running title:
Invariant faces associated with a cone-preserving map
\newpage
\baselineskip 0.6 cm

\vskip 60 pt
\begin{center}
{\large\bf ON THE INVARIANT FACES ASSOCIATED WITH A CONE-PRESERVING MAP}
\end{center}
\par
\vskip 12 pt
\begin{center}
{\small BIT-SHUN TAM
AND HANS SCHNEIDER}
\end{center}
\par\parindent=32 pt
\vskip 12 pt

\begin{minipage}[t]{11.5 cm}
{\small A}{\scriptsize BSTRACT.}$\,\,$
{\small For an $n\!\times\! n$ nonnegative matrix $P$,
an isomorphism is obtained between the lattice of initial subsets
(of $\{ 1,\cdots,n\}$) for $P$ and the lattice of $P$-invariant faces of the
nonnegative orthant $\IR^{\,\,n}_{\,\,+}$.
Motivated by this isomorphism,
we generalize some of the known combinatorial spectral results on a nonnegative
matrix that are given in terms of its classes to results for a cone-preserving
map on a polyhedral cone,
formulated in terms of its invariant faces.
In particular,
we obtain the following extension of the famous Rothblum Index Theorem for a
nonnegative matrix:
If $A$ leaves invariant a polyhedral cone $K$,
then for each distinguished eigenvalue $\lambda$ of $A$ for $K$,
there is a chain of
$m_\lambda$ distinct $A$-invariant join-irreducible faces of $K$,
each containing in its relative interior a generalized eigenvector of $A$
corresponding to $\lambda$ (referred to as
semi-distinguished $A$-invariant faces associated with $\lambda$),
where $m_\lambda$ is the maximal order of distinguished generalized
eigenvectors of $A$ corresponding to $\lambda$,
but there is no such chain with more than $m_\lambda$ members.
We introduce the important new concepts of semi-distinguished $A$-invariant
faces,
and of spectral pairs of faces associated with a cone-preserving map,
and obtain several properties of a cone-preserving map that
mostly involve these two concepts,
when the underlying cone is polyhedral,
perfect,
or strictly convex and/or smooth,
or is the cone of all real polynomials of degree not exceeding $n$ that are
nonnegative on a closed interval.
Plentiful illustrative examples are provided.
Some open problems are posed at the end.}
\end{minipage}
\par\parindent=0 pt
{
\unitlength 1mm
\begin{picture}(90,5)
\put(0,0){\line(1,0){40}}
\end{picture}
\par\parindent=16 pt
{\scriptsize {{1991 {\it Mathematics Subject Classification}.
Primary 15A48;
Secondary 47B65,
47A25,
46B42.
\par\parindent=16 pt
{\it Key words and phrases}.
Cone-preserving map,
nonnegative matrix,
polyhedral cone,
perfect cone,
strictly convex smooth cone,
spectral pair of a vector,
spectral pair of a face,
Perron-Schaefer condition,\par\parindent=0 pt
initial subset,
semi-distinguished class,
semi-distinguished invariant face,
distinguished generalized eigenvector,
chain of invariant faces.
\par\parindent=16 pt
\vskip -1 pt
Research of the first author partially supported by the National Science
Council of
the
Republic of \par\parindent=0 pt\vskip -2 pt China grant NSC 86-2115-M-032-002-;
the second author's research partially supported by NSF grants \par\vskip -3 pt
DMS-9123318 and
DMS-9424346 }}}
\par
\vfill{ }
\newpage

\par\parindent=0 pt
\begin{center}
1. I{\scriptsize NTRODUCTION}
\end{center}
\par\parindent=16 pt

This is the fourth of a sequence of papers in an attempt to study the
Perron-Frobenius theory of a (square,
entrywise) nonnegative matrix and its\par\parindent=0 pt
generalizations in the finite-dimensional setting from the cone-theoretic\par
viewpoint.
The first three papers in this sequence are [T--W],
[Tam 3] and \par [T--S 1].
Different fundamental aspects of a cone-preserving map are \par considered,
one in each paper of the sequence --- the Collatz-Wielandt sets (or numbers)
in [T--W],
the distinguished eigenvalues in [Tam 3],
the core in [T--S 1],
and the invariant faces in this paper.
\par\parindent=16 pt

In [T--S 1],
the preceding paper of our sequence,
we studied the core of a cone-preserving map $A$,
i.e. the convex cone $\core_K(A):=\bigcap^{\infty}_{i=1}A^iK$.
We were able to characterize the $K$-irreducibility or $K$-primitivity of $A$
in terms of its core,
and also to obtain many interesting results that give connections between
the core,
the peripheral spectrum,
the Perron-Schaefer condition and the distinguished $A$-invariant faces of $K$.
(We shall define some of these terms later on.)
Nevertheless,
the results we obtained directly or indirectly suggest that the set
$\core_K(A)$
does not capture all the important information about the spectral properties of
$A$.
In particular,
if $\core_K(A)$ is a polyhedral cone,
then it does not contain any distinguished generalized eigenvectors of $A$
other
than eigenvectors.
Thus,
the index of the spectral radius of $A$ cannot be determined from a knowledge
of its core.
On the other hand,
in the nonnegative matrix case,
as is known,
the index of the spectral radius or of a distinguished eigenvalue can be
described
in terms of its (equivalence) classes (of communicating states).
Moreover,
we found,
in a preliminary study,
a close connection between the classes of a nonnegative
matrix and its associated invariant faces (of the corresponding nonnegative
orthant).
This motivates our present study of the invariant faces associated with a
cone-preserving
map.
\par\parindent=16 pt

Here we would like to stress the differences between our work and that of
others on positive
(linear)
operators.
So far,
the deepest results in the theory of positive operators are mostly obtained
under the assumption that the positive operator under consideration is
irreducible or the underlying space is a nice space,
such as a Lebesgue space or,
more generally,
a \par\parindent=0 pt Banach lattice.
Furthermore,
many of these deep results have no counter-parts in the finite-dimensional
setting.
A look at the monographs [MN 2],
[Scha] or the recent survey paper by Dodds [Dod] will give the reader some
ideas.
(See also [Zer] for an earlier survey.)
Over the past two decades,
progress in a large part of the theory has followed three principal themes,
namely,
the kernel representation theorem of Bukhvalov (1974),
the compact majorization theorem of Dodds-Fremlin (1979) and the spectral
radius \par theorem of de Pagter (1986).
These results are given in the settings of the \par familiar Lebesgue spaces or
of a
Banach lattice and,
as can be readily seen none of them have a finite-dimensional counter-part,
as their conclusions
(that a certain operator is compact or is a kernel operator)
are trivially \par satisfied when the underlying space is of finite dimension.
Nevertheless,
the recent development of the combinatorial spectral theory of nonnegative
\par
matrices (applied to the reducible case) has had
some impact on the theory of positive operators.
For instance,
many of the graph-theoretic ideas
(such as the concepts of classes,
the notion of accessibility between states or classes)
or of the combinatorial spectral results
(such as the Nonnegative Basis theorem and the Rothblum Index theorem)
of a nonnegative matrix,
implicit in the early work of Schneider
([Schn 1])
and formulated first by Rothblum [Rot],
have already been extended,
first to the setting of an eventually compact linear integral operator on
$L^p(\Omega,\mu)$,
$1\le p<\infty$,
with nonnegative kernel
(see [Nel],
[Vic 1, 2],
[J--V 1]),
and then to the setting of a positive,
eventually compact linear operator on a Banach lattice having order continuous
norm
(see [J--V 2, 3]).
Their treatment is made possible by a decomposition of the positive operator
under consideration in terms of certain closed ideals of the underlying space
in a form which directly generalizes the Frobenius normal form of a nonnegative
reducible matrix.
In this work,
we are also generalizing theorems on reducible nonnegative matrices.
However,
we shall not generalize the Frobenius normal form;
indeed,
we do not expect there is a natural generalization for a cone-preserving map
when the underlying cone is a general proper cone.
Our focus is on the invariant faces,
and our results are geometric and provide new information in finite dimension.
Moreover,
we do not confine ourselves to a lattice-ordered cone,
which in the finite-dimensional setting reduces to the classical nonnegative
matrix case.
\par\parindent=16 pt

We now describe the contents of this paper in some detail.
Some necessary definitions and preliminary results are collected in Section 2.
\par\parindent=16 pt

For a nonnegative matrix,
we have the known concepts of basic classes,
\par\parindent=0 pt
distinguished classes,
initial classes,
etc.
In order to find for a cone-preserving map the appropriate analogs of these
concepts,
in Section 3 we begin our investigation by studying the
lattice of $P$-invariant faces
of $\IR_{\,\,+}^{\,\,n}$ (the \par nonnegative orthant of $\IR^{\,\,n}$)
associated with an $n\!\times\! n$ nonnegative matrix $P$.
We introduce the combinatorial concept of an initial subset (of
$\{ 1,\cdots,n\}$) for $P$.
An isomorphism is found between the lattice of initial subsets for $P$ and the
lattice of $P$-invariant faces of $\IR^{\,\,n}_{\,\,+}$,
thus relating a combinatorial concept to a geometric concept via algebraic
structures.
In particular,
under this isomorphism a $P$-invariant join-irreducible face which contains
in its
relative interior a generalized eigenvector of $P$ corresponding to $\rho(P)$,
the spectral radius of $P$,
corresponds to an initial subset determined by a basic class.
This isomorphism motivates the subsequent attempts in this paper to generalize
many of the known combinatorial spectral results on a nonnegative matrix that
are
given in terms of its classes to results for a cone-preserving map,
formulated in terms of its invariant faces.
In the course of establishing this isomorphism,
we reprove the well-known Frobenius$\,$-Victory theorem,
introduce the concept of a semi-distinguished class for a nonnegative
matrix,
reformulate the known nonnegative basis theorem for the generalized eigen-space
corresponding to a distinguished eigenvalue of a nonnegative matrix and offer a
cone-theoretic proof.
\par\parindent=16 pt

The following are two of the well-known results on the combinatorial
\par\parindent=0 pt
spectral theory
of a nonnegative matrix that we wish to extend in this paper.
\par\parindent=0 pt
\vskip 12 pt

{\bf Theorem A.}
(Rothblum Index Theorem) [Rot, Theorem 3.1(2)]\quad
{\it If $P$ is a nonnegative matrix,
then there always exists a chain of $\nu$ basic classes of $P$,
where $\nu$ is the index of $\rho(P)$ as an eigenvalue of $P$,
but there is no chain with more than $\nu$ basic classes.}
\par\parindent=0 pt
\vskip 12 pt

{\bf Theorem B.}
([Schn 1, Theorem 5] or [Schn 3, Theorem 6.3])\quad
{\it If $P$ is a nonnegative matrix,
then the Jordan canonical form of $P$ contains only one block corresponding to
$\rho(P)$
(i.e. $\dim{\cal N}(\rho(P)I-P)=1$)
if and only if any two basic classes of $P$ are comparable
(with respect to the accessibility relation).}
\par\parindent=16 pt
\vskip 12 pt

In our investigation,
we find it useful to introduce the concept of the \par\parindent=0 pt spectral
pair of a face relative
to a cone-preserving map.
We do this in \par Section 4.
Given an $n\!\times\! n$ complex matrix $A$,
for any vector $x\in {\IC}^n$,
we \par assign to it an ordered pair
$(\rho_{{}_{\scriptstyle x}}(A), \ord_A(x))$,
where $\rho_{{}_{\scriptstyle x}}(A)$ is the local spectral radius of $A$
at $x$ and $\ord_A(x)$
is the order of $x$ relative to $A$,
a new concept which we introduce as an extension of the usual concept of the
order
of a generalized eigenvector.
We denote this ordered pair by $\sp_A(x)$ and refer to it as the spectral pair
of $x$
relative to $A$.
If $A$ leaves invariant a proper cone $K$ and $F$ is a face of $K$,
we observe that the spectral pair $\sp_A(x)$ is \par\parindent=0 pt independent
of the choice of
$x$
from the relative interior of $F$.
This constant value is denoted by $\sp_A(F)$,
and will be referred to as the spectral pair of $F$ relative to $A$.
Fundamental properties of spectral pairs of faces
(or of vectors in the underlying cone)
relative to a cone-preserving map are derived.
One consequence is that
we obtain a class of $A$-invariant faces of $K$ given in terms of spectral
pairs
and the lexicographic ordering defined between ordered pairs of real numbers,
which extends and refines in the finite-dimensional case the class of invariant
ideals discovered by Meyer-Nieberg [MN 1] for a positive linear operator
defined
on a Banach lattice.
We also prove several equivalent conditions for a distinguished $A$-invariant
face,
a concept which we \par introduced in our preceding paper [T--S 1] as an
extension
of the concept
of a distinguished class of a nonnegative matrix.
Then as a natural extension of the concept of a semi-distinguished class we
call
an $A$-invariant face semi-distinguished if it is $A$-invariant
join-irreducible
and contains in its relative interior a generalized eigenvector of $A$.
\par\parindent=16 pt

In Section 5,
in terms of semi-distinguished $A$-invariant faces we \par\parindent=0 pt
formulate and establish
a natural extension of the Rothblum Index theorem in the
setting of a linear
mapping preserving a polyhedral cone.
We also give examples which show that in the nonpolyhedral case the result
is no
longer true.
\par\parindent=16 pt

We take a closer look at semi-distinguished $A$-invariant faces in Section 6.
Further properties of invariant faces
(mostly involving spectral pairs,
semi-distinguished invariant faces,
or distinguished generalized eigenvectors)
associated with a linear mapping preserving a polyhedral cone are obtained.
It is proved that if $A$ leaves invariant a polyhedral cone $K$,
then a nonzero $A$-invariant face $F$ of $K$ is semi-distinguished if and only
if
for any $A$-invariant face $G$ properly included in $F$,
we have,
$\sp_A(G)\prec\sp_A(F)$,
where $\prec$ denotes the strict lexicographic ordering between ordered pairs
of real numbers.
The nonpolyhedral case is also examined in detail.
It appears that whether a cone-preserving map shares some of the properties of
a cone-preserving
map on a polyhedral cone depends much on the geometry of the underlying cone.
In particular,
we prove that if $A$ leaves invariant a proper $K$ which has the property
that the dual cone of each of its faces is a facially exposed cone,
then the above-mentioned characterization of semi-distinguished
$A$-invariant faces in terms of lexicographic ordering still holds.
\par\parindent=16 pt

In Section 7 we obtain an extension and a refinement of Theorem B in the
setting
of a linear mapping preserving a polyhedral cone.
The situation for the nonpolyhedral case is also investigated carefully.
\par\parindent=16 pt

We give some multi-purpose examples in Section 8,
and pose some open questions in Section 9,
the final section of our paper.
\par\parindent=16 pt

Applications of our theory to solvability of linear equations over cones will
be given in a forthcoming paper [T--S 2].
\par\parindent=16 pt

Before we end this section,
we would like to point out that
in this paper our proofs of theorems on a linear mapping preserving a
polyhedral
cone rely much on the fact that if $A$ leaves invariant a polyhedral cone $K$,
then $K$ contains generalized eigenvectors of $A$ of all possible orders
between
$1$ and $\nu$,
where $\nu$ is the index of $\rho(A)$ as an eigenvalue of $A$.
The latter fact follows from [Tam 3, Theorem 7.5(ii)] which asserts that there
always
exists a vector $x$ in $K$
(a polyhedral cone) such that ${(A-\rho(A)I)}^\nu x={\bf 0}$,
and ${(A-\rho(A)I)}^ix$ are nonzero vectors of $K$ for $i=0,\cdots,\nu-1$.
This last fact is,
in turn,
proved in [Tam 3],
using the minimal generating matrix of a polyhedral cone as a tool and
depending
on the corresponding known result for a nonnegative matrix.
However,
in our next paper [Tam 4],
we shall offer an analytic and cone-theoretic proof of a result stronger than
this last fact,
without assuming the corresponding result for a nonnegative matrix.
Thus,
our treatment can be kept independent of the known results on a nonnegative
matrix.
\par\parindent=0 pt
\vskip 25 pt

\begin{center}
2. P{\scriptsize RELIMINARIES}
\end{center}
\par\parindent=16 pt

We shall restrict our attention to finite-dimensional vector spaces.
A familiarity with convex cones,
cone-preserving maps,
and graph-theoretic properties of nonnegative matrices is assumed.
(For references,
see [Bar 2],
[B--P] or [Schn 3].)
The terminology and notation of the previous papers in this sequence are
adopted here.
For convenience,
we collect some of them below,
together with some additional new definitions and notation.
\par\parindent=16 pt

Unless specified otherwise,
all matrices considered in this paper are square.
We always use $K$ to denote a proper
(i.e. closed,
full,
pointed convex)
cone of $\IR^{\,\,n}$,
and $\pi(K)$ to denote the set of all $n\!\times\! n$ real matrices $A$ that
satisfy $AK\subseteq K$.
Elements of $\pi(K)$ are referred to as
{\it cone-preserving maps}
(or more commonly as {\it positive-operators} on $K$).
\par\parindent=16 pt

If $S\subseteq K$,
we denote by $\Phi(S)$ the {\it face of $K$ generated by $S$},
that is,
the intersection of all faces of $K$ including $S$.
If $x\in K$,
we write $\Phi(\{ x\})$ simply as $\Phi(x)$.
A vector $x\in K$ is called an {\it extreme vector} if either $x$ is the zero
vector or $x$ is nonzero and $\Phi(x)=\{\lambda x: \lambda \ge 0\}$;
in the latter case,
the face $\Phi(x)$ is called an {\it extreme ray}.
A proper cone $K$ is said to be {\it polyhedral} if it has finitely many
extreme
rays.
\par\parindent=16 pt

The set of all faces of $K$ is denoted by ${\cal F}(K)$.
By a {\it nontrivial face} of $K$ we mean a face other than $K$ itself
or the zero face $\{ {\bf 0}\}$.
By the {\it dual cone} $K^*$ of $K$ we mean the proper cone
$\{ z\in\IR^{\,\,n}:\langle z,x\rangle\ge 0\,\,\mbox{for all}\,\, x\in K\}$,
where $\langle z,x\rangle$ is the usual inner product of $\IR^{\,\,n}$ between
the vectors $z$ and $x$.
We also use $d_K$ to denote the {\it duality operator of $K$},
i.e.
the mapping from ${\cal F}(K)$ to ${\cal F}(K^*)$ given by:
$d_K(F)={(\span F)}^\perp\bigcap K^*$
(see [Tam 2]).
A face $F$ of $K$ is said to be {\it exposed} if there exists a face $G$ of
$K^*$ such that
$F=d_{K^*}(G)$;
or equivalently,
if $F=d_{K^*}\circ d_K(F)$.
A cone $K$ is called {\it facially exposed} if each of its
faces is exposed;
or equivalently,
$d_K$ is injective,
or $d_{K^*}$ is surjective
(see [Tam 2,
Proposition 2.5 (a) and Corollary 2.6]).
\par\parindent=16 pt

The set of all $n\!\times\! n$ complex matrices is denoted by ${\cal M}_n$.
The range \par\parindent=0 pt space,
nullspace and the spectral radius of a matrix $A\in {\cal M}_n$ are denoted\par
respectively by ${\cal R}(A)$,
${\cal N}(A)$ and $\rho(A)$.
If $A$ is a real matrix,
its range space \par and null space are understood to be taken in the
corresponding real spaces.\par
Eigenvalues of $A$ with modulus $\rho(A)$ are said to compose the
{\it peripheral \par spectrum of $A$}.
\par\parindent=16 pt

For any $A\in {\cal M}_n$ and $b\in{\IC}^n$,
by the {\it cyclic space}
(relative to $A$)
{\it generated by the vector $b$},
denoted by $W_b$,
we mean the subspace $\span \{ b,Ab,A^2b,\cdots\}$.
(Again,
when $A$ and $b$ are both real,
the above linear span is taken in the real space.)
If $\lambda$ is an eigenvalue of $A$,
we denote by $\nu_{\lambda}(A)$ the index of $\lambda$ as an eigenvalue of $A$.
\par\parindent=16 pt

Let $A$ be an $n\!\times\! n$ real matrix.
It is known that a necessary and sufficient condition for the existence of a
proper
cone $K$ of $\IR^{\,\,n}$ such that $A\in\pi(K)$ is that,
for each eigenvalue $\lambda$ in the peripheral spectrum of $A$,
we have $\nu_{\lambda}(A)\le \nu_{\rho(A)}(A)$.
(Then clearly $\rho(A)$ is an eigenvalue of $A$.)
This condition is now referred to as the {\it Perron-Schaefer condition}
(see [Schn 2, the paragraph after Theorem 1.1]).
\par\parindent=16 pt

If $A\in {\cal M}_n$ and $x$ is a nonzero vector of $\IC^n$ such that
${(\lambda I-A)}^kx={\bf 0}$ for some $\lambda\in\IC$ and some positive integer
$k$,
then $x$ is called a {\it generalized eigenvector of $A$} corresponding to
(the eigenvalue) $\lambda$.
The least such integer $k$ is called the {\it order of $x$} as a generalized
eigenvector of $A$,
and is denoted by $\ord_A(x)$.
\par\parindent=16 pt

For any $A\in{\cal M}_n$ and $x\in{\IC}^n$,
we define the {\it local spectral radius of $A$ at $x$},
denoted by $\rho_{{}_{\scriptstyle x}}(A)$,
as follows.
If $x$ is the zero vector,
take $\rho_{{}_{\scriptstyle x}}(A)$ to be $0$.
If $x$ is a nonzero vector,
write $x$ uniquely as a sum of generalized eigenvectors of $A$,
say,
\begin{center}
$x=x_1+\cdots+x_k$,
\end{center}

\par\parindent=0 pt
where $k\ge 1$,
and $x_1,\cdots,x_k$ are generalized eigenvectors of $A$ corresponding to
distinct
eigenvalues $\lambda_1,\cdots,\lambda_k$ respectively.
Define $\rho_{{}_{\scriptstyle x}}(A)$ to be $\max_{1\le i\le
k}\left|\lambda_i\right|$.
For other equivalent definitions,
see [T--W, Theorem 2.2].
\par\parindent=16 pt

If $A\in\pi(K)$ and $x\in K$ is an eigenvector
(respectively,
generalized \par\parindent=0 pt eigenvector),
then $x$ is called a {\it distinguished eigenvector}
(respectively,\par
{\it distinguished generalized eigenvector})
{\it of $A$ for $K$},
and the corresponding eigenvalue is known as a
{\it distinguished eigenvalue of $A$ for $K$}.
When there is no danger of confusion,
we simply use the terms distinguished eigenvector,
distinguished generalized eigenvector and distinguished eigenvalue (of $A$).
(For reference,
see [Tam 3].)
By an {\it extremal distinguished eigenvector}
of $A$ corresponding to the (distinguished) eigenvalue $\lambda$,
we mean a nonzero extreme vector of the cone ${\cal N}(\lambda I-A)\bigcap K$.
\par\parindent=16 pt

The nonnegative orthant of $\IR^{\,\,n}$ is denoted by $\IR^{\,\,n}_{\,\,+}$.
It is clear that $\IR^{\,\,n}_{\,\,+}$ is a proper cone of $\IR^{\,\,n}$,
and also that $\pi(\IR^{\,\,n}_{\,\,+})$ is the cone of all $n\!\times\! n$
nonnegative matrices.
When we consider an $n\!\times\! n$ nonnegative matrix as a cone-preserving
map,
unless specified otherwise,
we always mean that the underlying cone is $\IR^{\,\,n}_{\,\,+}$.
\par\parindent=16 pt

A familiarity with the diagraph and classes of a nonnegative matrix is assumed.
(For reference,
see,
for instance,
[Schn 3].)
If $n$ is a positive integer,
we denote by $\langle n\rangle$ the set $\{ 1,2,\cdots,n\}$.
We always use $P$ to denote an $n\!\times\! n$ nonnegative matrix for some
positive integer $n$.
The subsets of $\langle n\rangle$ are usually denoted by capital Latin letters
$I$,
$S$,
etc.,
except that we use small Greek letters $\alpha$,
$\beta$,
etc.\ to denote the classes of $P$.
The accessibility relation is usually defined between the classes of $P$,
but we also say $i$ {\it has access to} $j$
(where $i, j\in\langle n\rangle$)
if there is a directed path in $G(P)$
(the {\it digraph of $P$}) from $i$ to $j$.
Similarly,
we also say $i$ has access to a subset of $\langle n\rangle$ with the obvious
meaning.
If $S, T\subseteq\langle n\rangle$,
we denote by $P_{ST}$ the principal submatrix of $P$ with rows indexed by $S$
and columns indexed by $T$.
\par\parindent=16 pt

Sometimes we say a class $\alpha$ of $P$ is associated with $\lambda$.
By that we mean $\lambda=\rho(P_{\alpha\alpha})$.
A class $\alpha$ is said to be {\it semi-distinguished} if
$\rho(P_{\beta\beta})\le\rho(P_{\alpha\alpha})$ for any class $\beta$ which
has access to $\alpha$ but not equal to $\alpha$.
If,
in this definition,
the weak inequality is replaced by strict inequality of the same type,
we recover the usual definition of a distinguished class.
\par\parindent=16 pt

A subset $I$ of $\langle n\rangle$ is called an {\it initial subset for $P$}
if either $I$ is empty,
or $I$ is nonempty and $P_{I^\prime I}={\bf 0}$,
where $I^\prime=\langle n\rangle\backslash I$;
equivalently,
for every $j\in\langle n\rangle$,
$I$ contains $j$ whenever $j$ has access to $I$.
Clearly an initial subset for $P$ which is also a class of $P$ is simply an
initial
class in the usual sense.
It is also not difficult to see that a nonempty subset $I$ of $\langle
n\rangle$
is an initial subset for $P$ if and only if $I$ is the union of an
{\it initial collection of classes} of $P$,
where a nonempty collection of classes of $P$ is said to be initial if whenever
it contains a class $\alpha$,
it also contains all classes having access to $\alpha$.
One can readily verify that the intersection and union of two initial subsets
for $P$ are both initial subsets for $P$.
Thus the collection of all initial subsets for $P$ forms a sublattice of
$2^{\langle n\rangle}$,
the lattice of all subsets of $\langle n\rangle$.
We denote this sublattice by $\cal I$.
\par\parindent=16 pt

It is known that ${\cal F}(K)$
forms a lattice under inclusion as the partial \par\parindent=0 pt ordering,
with meet and join given by:
$F\bigwedge G=F\bigcap G$ and $F\bigvee G=\Phi(F\bigcup G)$.
If $A\in\pi(K)$,
then a face $F$ of $K$ is said to be {\it $A$-invariant} if $AF\subseteq F$.
Clearly,
the meet of two $A$-invariant faces is $A$-invariant.
The join of two $A$-invariant faces is also $A$-invariant;
this is because,
for any two such faces $F$,
$G$,
we have
\begin{center}
$A(F\bigvee G)=A(\Phi(F\bigcup G))\subseteq\Phi(A(F\bigcup
G))\subseteq\Phi(F\bigcup G)=F\bigvee G$.
\end{center}
\par\parindent=0 pt
Thus,
the set of all $A$-invariant faces of $K$,
which we denote by ${\cal F}_A(K)$
(or simply by ${\cal F}_A$ if there is no danger of confusion),
forms a sublattice of ${\cal F}(K)$.
A face $F$ of $K$ is said to be {\it $A$-invariant join-reducible} if $F$ is
$A$-invariant and is the join of two (and hence,
of all) $A$-invariant faces of $K$ that are properly included in $F$;
in other words,
$F$ is join-reducible in the lattice ${\cal F}_A$ in the usual
lattice-theoretic sense.
(For definitions of lattice-theoretic terms,
see [Bir].)
An $A$-invariant face which is not $A$-invariant join-reducible is said to be
{\it $A$-invariant join-irreducible}.
According to our definition,
the zero face is always $A$-invariant join-irreducible.
\par\parindent=16 pt

For any $F\in{\cal F}_A$,
we denote by $\rho_{{}_{\scriptstyle F}}(A)$ (or simply by
$\rho_{{}_{\scriptstyle F}}$) the spectral radius of the
restriction map ${\left. A\right|}_{\span F}$;
then we also say $F$ is {\it associated with} $\rho_{{}_{\scriptstyle F}}$.
A face $F$ is said to be a {\it distinguished $A$-invariant face of $K$}
({\it associated with} $\lambda$) if $F$ is a nonzero $A$-invariant face of $K$
such that $\rho_{{}_{\scriptstyle G}}<\rho_{{}_{\scriptstyle F}}$ for any
nonzero $A$-invariant face $G$ properly included
in $F$ (and $\rho_{{}_{\scriptstyle F}}=\lambda$).
\par\parindent=16 pt

We shall need the following explicit characterization of the smallest
$A$-invariant
face that contains a given vector $x$ of $K$.
\par\parindent=0 pt
\vskip 12 pt

{\bf Lemma 2.1.}
{\it Let $A\in\pi(K)$ and let $x\in K$.
Then $\Phi({(I+A)}^{n-1}x)$ is the smallest $A$-invariant face of $K$
containing $x$.}
\par\parindent=0 pt
\vskip 12 pt

{\it Proof.}
For convenience,
denote the vector ${(I+A)}^{n-1}x$ by $y$.
Note that $y$ is a positive linear combination of the vectors
$x,Ax,\cdots,A^{n-1}x$
of $K$ and also that these vectors span $W_x$
(the cyclic space generated by $x$).
So $y$ belongs to the relative interior of $W_x\bigcap K$.
But clearly $Ay\in W_x\bigcap K$,
so there exists some $\alpha>0$ such that $\alpha y-Ay\in W_x\bigcap K$,
and hence $\alpha y\ge^K Ay$.
Thus,
$\Phi(y)$ is an $A$-invariant face which,
in addition,
also contains $x$,
as $y\ge^K x$.
If $F$ is an $A$-invariant face of $K$ that contains $x$,
then $F$ must contain $Ax,\cdots,A^{n-1}x$,
and hence also $y$;
thus $F\supseteq\Phi(y)$.
This proves that $\Phi(y)$ is the smallest $A$-invariant face of $K$ that
contains $x$.
\par
\hfill{$\Box$}
\par\parindent=0 pt
\vskip 25 pt
\newpage

\begin{center}
3. T{\scriptsize HE INVARIANT FACES

ASSOCIATED WITH A NONNEGATIVE MATRIX}
\end{center}
\par\parindent=16 pt

It is known that every face of $\IR^{\,\,n}_{\,\,+}$ is of the form
\begin{center}
$F_I=\{ x\in\IR^{\,\,n}_{\,\,+}: \supp(x)\subseteq I\}$,
\end{center}
\par\parindent=0 pt
where $I\subseteq\langle n\rangle$,
and $\supp(x)$ is the {\it support of $x$},
i.e. the set $\{ i\in\langle n\rangle: \xi_i\not= 0\}$ for
$x={(\xi_1,\cdots,\xi_n)}^T$.
Indeed,
as can be readily seen,
the mapping $I\longmapsto F_I$ is an isomorphism between the lattices
$2^{\langle n\rangle}$ and ${\cal F}(\IR^{\,\,n}_{\,\,+})$.
The following result asserts that if $P$ is an $n\!\times\! n$ nonnegative
matrix,
then under this isomorphism,
the initial subsets for $P$ correspond to the $P$-invariant faces of
$\IR^{\,\,n}_{\,\,+}$.
\par\parindent=0 pt
\vskip 12 pt

{\bf Theorem 3.1.}
{\it Let $P$ be an $n\!\times\! n$ nonnegative matrix.
Denote by $\cal I$ the lattice of all initial subsets for $P$
and by
${\cal F}_P$ the lattice of all $P$-invariant faces of $\IR^{\,\,n}_{\,\,+}$.
Then the association $I\longmapsto F_I$ induces an isomorphism from the lattice
$\cal I$
onto the lattice ${\cal F}_P$.}
\par\parindent=0 pt
\vskip 12 pt

{\it Proof.}
Let $\triangle: {\cal I}\longrightarrow{\cal F}_P$ be the mapping defined by
$\triangle(I)=F_I$.
First, we need to show that $\triangle$ is a well-defined mapping.
Let $I\in{\cal I}$.
To show that $F_I$ is $P$-invariant,
take any vector $x\in F_I$.
Then
\begin{tabbing}
\qquad\qquad\quad$\supp(Px)$\= $=\{ j\in\langle n\rangle: p_{ji}\not= 0$ for
some $i\in\supp(x)\}$\\
\> $\subseteq\{ j\in\langle n\rangle: j$ has access to $\supp(x)\}$\\
\> $\subseteq I$,
\end{tabbing}
\par\parindent=0 pt
where the last inclusion follows from the fact that $\supp(x)\subseteq I$
and $I$ is an initial subset for $P$.
Hence,
by definition,
we have $Px\in F_I$.
\par\parindent=16 pt

Because $\triangle$ is just a restriction of the isomorphism $I\longmapsto F_I$
between the lattices $2^{\langle n\rangle}$ and ${\cal
F}(\IR^{\,\,n}_{\,\,+})$,
and $\cal I$ and ${\cal F}_P$ are respectively sublattices of
$2^{\langle n\rangle}$ and ${\cal F}(\IR^{\,\,n}_{\,\,+})$,
$\triangle$ is a one-to-one lattice homomorphism.
It remains to show that the mapping $\triangle$ is onto.
\par\parindent=16 pt

Let $F$ be a $P$-invariant face.
As a face of $\IR^{\,\,n}_{\,\,+}$,
$F$ must be of the form $F_S$ for some $S\subseteq\langle n\rangle$;
indeed,
$S$ is the set $\bigcup_{x\in F}\supp(x)$.
Let $i\in S$.
Consider any $j\in\langle n\rangle$,
which has access to $i$ but not equal to $i$.
Since $j$ has access to $i$,
there exists a positive integer $k$ such that the $(j,i)$ entry of $P^k$ is
nonzero;
or,
in other words,
$j\in\supp(P^ke_i)$,
where $e_i$ is the ith standard unit vector of $\IR^{\,\,n}$.
But $P^ke_i\in F$,
as $e_i\in F$ and $F$ is $P$-invariant,
so $j\in S$.
This shows that $S$ is an initial subset for $P$,
and thus establishes the surjectivity of $\triangle$.
The proof is complete.
\par
\hfill{$\Box$}
\par\parindent=16 pt
\vskip 12 pt

Let $\cal K$ be a nonempty collection of classes of $P$.
A class $\alpha\in{\cal K}$ is said to be {\it final in} $\cal K$ if it has no
access to
other classes of $\cal K$.
If $I$ is a nonempty initial subset for $P$,
then as mentioned before,
$I$ is the union of all classes in an initial collection $\cal K$ of classes of
$P$.
An initial subset $I$ for $P$ is said to be {\it determined by a class}
$\alpha$ of $P$ if it is the union of all classes of $P$ having access to
$\alpha$
(or,
in other words,
$\alpha$ is the only class final in the initial collection corresponding to
$I$);
in this case,
we write $F_I$ simply as $F_\alpha$.
The use of the notation $F_\alpha$ here agrees with that of [T--S 1, Section
4].
Since we use only small Greek letters to denote classes,
the use of this notation should cause no confusion.
According to Theorem 3.1,
each $F_\alpha$ is a $P$-invariant face.
\par\parindent=0 pt
\vskip 12 pt

{\bf Lemma 3.2.}
{\it Let $P$ be an $n\!\times\! n$ nonnegative matrix.
If $\alpha$ is a distinguished class of $P$ associated with the eigenvalue
$\lambda$,
then $F_\alpha$ contains a unique (up to multiples) eigenvector of $P$
corresponding to $\lambda$.
Furthermore,
this eigenvector belongs to $\relint F_\alpha$ and is an extremal distinguished
eigenvector of $P$.}
\par\parindent=0 pt
\vskip 12 pt

{\it Proof.}
First,
note that ${\left. P\right|}_{\span F_\alpha}$ can be represented by the
principal
\par\parindent=0 pt submatrix of $P$ with rows and columns indexed by the
initial subset for $P$ determined
by $\alpha$.
As such,
the spectrum of ${\left. P\right|}_{\span F_\alpha}$ is simply the union of the
spectra of all possible submatrices $P_{\beta\beta}$ for which $\beta$ is a
class
of $P$ that has access to $\alpha$.
Since $\alpha$ is a distinguished class,
it follows that $\rho({\left. P\right|}_{\span
F_\alpha})=\rho(P_{\alpha\alpha})=\lambda$.
According to the Perron-Frobenius theorem,
there exists an eigenvector $x$ of $P$ corresponding to $\lambda$ that lies in
$F_\alpha$.
If $x\in\rbd F_\alpha$,
then $\Phi(x)$ is a $P$-invariant face properly included in $F_\alpha$,
and by Theorem 3.1 there is an initial subset $I$ for $P$ properly included in
the
initial subset determined by $\alpha$ such that $\Phi(x)=F_I$.
But then $\rho_{{}_{\scriptstyle {\Phi(x)}}}=\rho(P_{II})<\lambda$,
which contradicts our assumption on $x$.
This shows that $x\in\relint F_\alpha$.
\par\parindent=16 pt

If $P$ has an eigenvector in $F_\alpha$ corresponding to $\lambda$,
which is not a multiple of $x$,
then we readily obtain an eigenvector corresponding to $\lambda$ that lies in
$\rbd F_\alpha$.
But as we have just shown,
this would lead to a contradiction.
This proves the uniqueness of the desired eigenvector $x$.
It also follows that $x$ is an extremal distinguished eigenvector.
\par
\hfill{$\Box$}
\par\parindent=16 pt
\vskip 12 pt

We now offer a cone-theoretic proof of the following Frobenius$\,$-Victory
theorem,
first named in [T--S 1] (see Theorem 2.1 and the subsequent \par\parindent=0 pt
remarks there).
\par\parindent=0 pt
\vskip 12 pt

{\bf Theorem 3.3.}
{\it Let $P$ be an $n\!\times\! n$ nonnegative matrix.}
\par\parindent=16 pt
(i)
{\it For any real number $\lambda$,
$\lambda$ is a distinguished eigenvalue of $P$ (for $\IR^{\,\,n}_{\,\,+}$)
if and only if there exists a distinguished class $\alpha$ of $P$ such that
$\rho(P_{\alpha\alpha})=\lambda$.}
\par\parindent=16 pt
(ii)
{\it If $\alpha$ is a distinguished class of $P$,
then there is a (up to multiples) unique nonnegative eigenvector
$x^\alpha={(\xi_1,\cdots,\xi_n)}^T$ corresponding to $\rho(P_{\alpha\alpha})$
with the property that $\xi_i>0$ if and only if $i$ has access to $\alpha$.}
\par\parindent=16 pt
(iii)
{\it For each distinguished eigenvalue $\lambda$ of $P$,
the cone ${\cal N}(\lambda I-P)\bigcap\IR^{\,\,n}_{\,\,+}$ is simplicial and
its
extreme vectors are precisely all the distinguished \par\parindent=0 pt
eigenvectors of $P$ of the form
$x^\alpha$ as given in} (ii),
{\it where $\alpha$ is a distinguished class such that
$\rho(P_{\alpha\alpha})=\lambda$.}
\par\parindent=0 pt
\vskip 12 pt

{\it Proof.}
\par\parindent=16 pt
(i):
The $\lq\lq$ if $"$ part follows from Lemma 3.2.
To prove the $\lq\lq$ only if $"$ part,
choose a distinguished $A$-invariant face $F$ of $\IR^{\,\,n}_{\,\,+}$
associated with $\lambda$.
By \par\parindent=0 pt Theorem 3.1 $F=F_I$ for some nonempty initial subset $I$
for $P$.
Note that there can be only one class final in the initial collection of
classes
\par corresponding to $I$,
say $\alpha$,
(i.e. $F_I=F_\alpha$) and also that $\alpha$ is a distinguished class
associated
with $\lambda$;
otherwise,
it would be possible to choose a nonempty initial subset $I^\prime$ properly
included in $I$,
and hence $F_{I^\prime}$ properly included in $F_I$,
such that $\rho_{{}_{\scriptstyle {F_{I^\prime}}}}=\rho_{{}_{\scriptstyle
{F_I}}}=\lambda$,
contradicting the assumption that $F$ is a distinguished $A$-invariant face.
\par\parindent=16 pt
(ii):
Follows from Lemma 3.2 (and the definition of $F_\alpha$).
\par\parindent=16 pt
(iii):
By Lemma 3.2 each distinguished eigenvector of the form $x^\alpha$,
where $\alpha$ is a distinguished class associated with $\lambda$,
is an extreme vector of ${\cal N}(\lambda I-P)\bigcap\IR^{\,\,n}_{\,\,+}$.
Conversely,
if $x$ is an (nonzero) extreme vector of ${\cal N}(\lambda
I-P)\bigcap\IR^{\,\,n}_{\,\,+}$,
then as can be readily shown,
$\Phi(x)$ is a distinguished $A$-invariant face of $\IR^{\,\,n}_{\,\,+}$.
>From the proof of the $\lq\lq$ only if $"$ part of (i),
then we have $\Phi(x)=F_\alpha$ for some distinguished class $\alpha$
associated
with $\lambda$.
But by Lemma 3.2 the eigenvector of $P$ corresponding to $\lambda$ that
lies in $F_\alpha$ is unique,
so $x=x^\alpha$.
This proves that the extreme vectors of the cone ${\cal N}(\lambda
I-P)\bigcap\IR^{\,\,n}_{\,\,+}$
are precisely those of the given form.
That the cone ${\cal N}(\lambda I-P)\bigcap\IR^{\,\,n}_{\,\,+}$ is simplicial,
i.e. its set of distinct (up to multiples) extreme vectors is linearly
independent,
follows from the combinatorial properties of the supports of these vectors.
\par
\hfill{$\Box$}
\par\parindent=16 pt
\vskip 12 pt

Using the Frobenius$\,$-Victory theorem,
in [Tam 3, Theorem 5.2]
a cone-theoretic proof is offered for the well-known Nonnegative Basis theorem
for a nonnegative matrix (or a singular $M$-matrix).
[A stronger result known as the Preferred-Basis theorem has also been found.
See [H--S].]
Hershkowitz and Schneider [H--S, Section 5] have also extended the Nonnegative
Basis theorem to cover the case of a distinguished eigenvalue.
We reformulate their result as Theorem 3.4 below and offer an independent proof
which is based on the Nonnegative-Basis theorem and makes use of cone-theoretic
arguments.
\par\parindent=0 pt
\vskip 12 pt

{\bf Theorem 3.4.}
{\it Let $P$ be an $n\!\times\! n$ nonnegative matrix,
and let $\lambda$ be a \par\parindent=0 pt distinguished eigenvalue of $P$.
Let $\alpha_1,\cdots,\alpha_p$ be all the semi-distinguished classes of $P$
associated with $\lambda$.
Then the generalized eigenspace of $P$ \par corresponding to $\lambda$ contains
nonnegative
vectors $x^{\alpha_1},\cdots,x^{\alpha_p}$ such that
$x^{\alpha_i}\in\relint F_{\alpha_i}$ for each $i$.
Furthermore,
any such collection forms a basis for the subspace spanned by the nonnegative
generalized eigenvectors of $P$ \par corresponding to $\lambda$.}
\par\parindent=16 pt
\vskip 12 pt

Before we come to the proof of Theorem 3.4,
we need the following lemma.
\par\parindent=0 pt
\vskip 12 pt

{\bf Lemma 3.5.}
{\it Let $P$ be an $n\!\times\! n$ nonnegative matrix.
If $F$ is a $P$-invariant face of $\,\IR^{\,\,n}_{\,\,+}$ which contains in its
relative interior a generalized eigenvector of $P$,
then $F$ must be of the form $F_I$,
where $I$ is a nonempty initial subset for $P$ with the property that each
class
final in the initial collection of classes corresponding to $I$ is a
semi-distinguished
class associated with $\rho_{{}_{\scriptstyle F}}$.}
\par\parindent=0 pt
\vskip 12 pt

{\it Proof.}
By Theorem 3.1 $F$ is of the form $F_I$ for some initial subset $I$ for $P$.
Denote by $B$ the principal submatrix of $P$ with columns and rows indexed by
$I$.
Since $P$ has a generalized eigenvector that lies in $\relint F$,
the nonnegative matrix $B$ has a positive generalized eigenvector.
Hence,
by [Tam 3, Theorem 5.1] $\rho_{{}_{\scriptstyle F}} (=\rho(B))$ is the only
distinguished eigenvalue
of $B^T$ (for the \par\parindent=0 pt corresponding nonnegative orthant).
By the Frobenius$\,$-Victory theorem it follows that every initial class of
$B^T$,
and hence every final class of $B$,
is associated with $\rho_{{}_{\scriptstyle F}}$.
But the final classes of $B$ are precisely the classes final in the initial
collection of classes corresponding to $I$,
hence our assertion follows.
\par
\hfill{$\Box$}
\par\parindent=0 pt
\vskip 12 pt

{\it Proof of Theorem 3.4.}
Let $\alpha$ be any semi-distinguished  class of $P$ associated with $\lambda$.
Let $Q$ be the principal submatrix of $P$ with rows and columns indexed by the
initial subset determined by $\alpha$.
Then
\begin{tabbing}
\qquad\qquad\qquad\quad$\rho(Q)$\= $=\max\{ \rho(P_{\beta\beta}): \beta$ has
access to $\alpha\}$\\
\> $=\lambda$,
\end{tabbing}
\par\parindent=0 pt
where the last equality follows from the assumption that $\alpha$ is a
semi-\par\parindent=0 pt distinguished class of $P$.
An application of the first part of the Nonnegative-Basis theorem
(see [Tam 3, Theorem 5.2]) to $Q$ guarantees that there exists a positive
generalized eigenvector $w$ of $Q$ corresponding to $\lambda$.
For convenience,
we index the components of $w$ by the initial subset of $P$ determined by
$\alpha$.
Let $x^\alpha$ be the vector in $\IR^{\,\,n}_{\,\,+}$ defined by:
${(x^\alpha)}_i$ equals $w_i$ if $i$ has access to $\alpha$,
and equals zero otherwise.
Then $x^\alpha$ is the desired generalized eigenvector of $P$ corresponding to
$\lambda$ which is associated with the semi-distinguished class $\alpha$.
\par\parindent=16 pt

It is clear that any collection of vectors $x^\alpha$,
with one vector for each semi-distinguished class $\alpha$ associated with
$\lambda$,
forms a linearly independent family,
in view of the combinatorial properties of these vectors.
To complete the proof,
it remains to show that if $u$ is a nonnegative generalized eigenvector of $P$
corresponding to $\lambda$,
then $u$ is spanned by any such collection of vectors $x^\alpha$.
For that purpose,
it suffices to show that $\supp(u)\subseteq I$,
where we use $I$ to denote the union of all classes of $P$ having access to at
least one semi-distinguished class associated with $\lambda$.
This is because,
then we can apply the Nonnegative Basis theorem to the principal submatrix of
$P$ with rows and columns indexed by $I$ to draw the desired conclusion.
\par\parindent=16 pt

Let $u$ be any nonnegative generalized eigenvector of $P$ corresponding to
$\lambda$.
By Lemma 2.1 the face $\Phi({(I+P)}^{n-1}u)$ is the smallest $P$-invariant face
of $\IR^{\,\,n}_{\,\,+}$ containing $u$,
and by Theorem 3.1 it must be of the form $F_S$,
for a unique initial subset $S$ for $P$.
(In fact,
$S$ is the smallest initial subset for $P$ including $\supp(u)$,
i.e. it is the union of all classes of $P$ having access to $\supp(u)$.)
Since $F_S$ contains a generalized eigenvector of $P$ corresponding to
$\lambda$
in its relative interior,
namely,
the vector ${(I+P)}^{n-1}u$,
by Lemma 3.5 each class final in the initial collection of classes
corresponding
to $S$ is a semi-distinguished class associated with $\lambda$.
Hence,
we have
\begin{center}
$\supp(u)\subseteq\supp({(I+P)}^{n-1}u)=S\subseteq I$,
\end{center}
\par\parindent=0 pt
as desired.
\par
\hfill{$\Box$}
}

{\baselineskip 0.6 cm
\par\parindent=16 pt
\vskip 12 pt

The following characterizations of different types of $P$-invariant faces
\par\parindent=0 pt
associated with a nonnegative matrix $P$ motivate much of our subsequent work
in this paper.
\par\parindent=0 pt
\vskip 12 pt

{\bf Theorem 3.6.}
{\it Let $P$ be an $n\!\times\! n$ nonnegative matrix.
Let $I$ be an initial subset for $P$.
Then}
\par\parindent=16 pt
(i)
{\it $F_I$ is a minimal nonzero $P$-invariant face if and only if $I$ is an
initial
class of $P$.}
\par\parindent=16 pt
(ii)
{\it $F_I$ is a nonzero $P$-invariant join-irreducible face if and only if
$I$ is an initial subset determined by a single class.}
\par\parindent=16 pt
(iii)
{\it $F_I$ is a $P$-invariant face which contains in its relative interior
a \par\parindent=0 pt generalized eigenvector
(respectively,
an eigenvector)
of $P$ corresponding to $\lambda$ if and only if $I$ is a nonempty initial
subset
such that each class final in the initial collection of classes corresponding
to
$I$ is a semi-distinguished
(respectively,
distinguished)
class associated with $\lambda$.}
\par\parindent=16 pt
(iv)
{\it $F_I$ is a $P$-invariant join-irreducible face which contains in its
\par\parindent=0 pt relative
interior a generalized eigenvector
(respectively,
an eigenvector)
of $P$ corresponding to $\lambda$ if and only if $I$ is an initial subset
determined by a semi-distinguished
(respectively,
distinguished)
class associated with $\lambda$.}
\par\parindent=0 pt
\vskip 12 pt

{\it Proof.}
As in the proof of Theorem 3.1 we denote by $\triangle$ the isomorphism from
$\cal I$,
the lattice of all initial subsets for $P$,
onto ${\cal F}_P$,
the lattice of all $P$-invariant faces of $\IR^{\,\,n}_{\,\,+}$,
given by $\triangle(I)=F_I$.
\par\parindent=16 pt
(i):
The minimal nonempty members of $\cal I$ are precisely those initial
\par\parindent=0 pt subsets of $P$
determined by an initial class.
Under the lattice isomorphism $\triangle$,
they clearly correspond to the minimal nonzero $P$-invariant faces of
$\IR^{\,\,n}_{\,\,+}$.
\par\parindent=16 pt
(ii):
Note that the nonempty join-irreducible members of the lattice $\cal I$ are
precisely
those initial subsets determined by a single class,
whereas the nonzero join-irreducible members of ${\cal F}_P$ are the nonzero
$P$-invariant
join-irreducible faces.
But $\triangle$ sends the empty set to the zero face and as a lattice
isomorphism
it preserves the join-irreducibility property,
so our assertion follows.
\par\parindent=16 pt
(iii):
$\lq\lq$ Only if $"$ part:
The $\lq\lq$ generalized eigenvector $"$ case follows \par\parindent=0 pt
immediately from Lemma 3.5.
To prove the $\lq\lq$ eigenvector $"$ case,
we denote by $B$,
as in the proof of Lemma 3.5,
the principal submatrix of $P$ with columns and rows indexed by $I$.
Here,
instead of using [Tam 3, Theorem 5.1],
we apply [Tam 3, Theorem 5.4] to $B$ and deduce that
$\nu_{\rho_{{}_{\scriptscriptstyle F}} }(B)=1$,
and also that $I$ is a nonempty initial subset with the property that each
class final in the initial collection of classes corresponding to $I$ is a
semi-distinguished class (of $B$,
or equivalently,
of $P$) associated with $\rho_{{}_{\scriptstyle F}}$.
It remains to show that each class final in the initial collection
corresponding
to $I$ is,
in fact,
a distinguished class.
(Here,
if we apply the Rothblum Index theorem,
we can already draw the desired conclusion.
But,
providing an alternative cone-theoretic proof for the Rothblum Index theorem
and its extension is one of the purposes of this paper,
so we proceed instead as follows.)
Suppose that there is a class,
say $\beta$,
which is final in the initial collection corresponding to $I$ but is not a
distinguished class.
Then by Theorem 3.4 we can find a nonnegative generalized eigenvector
$x={(\xi_i)}^T$
of $B$ corresponding to $\rho_{{}_{\scriptstyle F}}$
with the property that $\xi_i>0$ if and only if
$i$ has access to $\beta$.
>From a consideration of supports of vectors,
it is clear that $x$ cannot be written as a linear combination of vectors of
the
form $x^\alpha$ as given in Theorem 3.3(ii),
where $\alpha$ is a distinguished class of $B$ associated with
$\rho_{{}_{\scriptstyle F}}$.
Hence,
it follows that $x$ is a generalized eigenvector of order greater than $1$.
Thus,
we obtain $\nu_{\rho_{{}_{\scriptscriptstyle F}} }(B)\ge 2$,
which is a contradiction.
\par\parindent=16 pt
$\lq\lq$ If $"$ part:
We give the proof for the $\lq\lq$ generalized eigenvector $"$ case,
and omit the proof for the $\lq\lq$ eigenvector $"$ case,
which depends on the use of the Frobenius$\,$-Victory theorem and is similar
and easier.
We denote by $\cal K$ the initial collection of classes of $P$ corresponding to
the
initial subset $I$.
\par\parindent=16 pt
Let $\alpha_1,\cdots,\alpha_k$ be the classes final in $\cal K$.
By our assumption on $I$,
each $\alpha_i$ is a semi-distinguished class associated with the eigenvalue
$\lambda$.
By Theorem 3.4,
for each $i$,
$1\le i\le k$,
there is a generalized eigenvector,
say $x_i$,
of $P$ corresponding to $\lambda$ that lies in relint $F_{\alpha_i}$.
Since $I$ is an initial subset,
it follows that the vector $\sum^k_{i=1}x_i$ is a generalized eigenvector of
$P$
corresponding to $\lambda$ that lies in $\relint F_I$.
\par\parindent=16 pt
(iv):
The $\lq\lq$ only if $"$ part and the $\lq\lq$ if $"$ part follow respectively
from the corresponding parts of (ii) and (iii).
\par
\hfill{$\Box$}
\par\parindent=0 pt
\vskip 25 pt

\begin{center}
4. T{\scriptsize HE SPECTRAL PAIR OF A FACE}
\end{center}
\par\parindent=16 pt

Let $A\in\pi(K)$.
For any face $F$ of $K$,
$A$-invariant or not,
we are going to assign to it an ordered pair of nonnegative real numbers,
known as the spectral pair of $F$ relative to $A$.
Such spectral pairs of faces will turn out to be useful in our later
investigations.
But first we need the concept of a spectral pair of a vector relative to an
arbitrary (square) matrix,
and for this purpose we introduce the concept of the order of a vector relative
to a matrix.
We thus extend the usual concept of the order of a generalized eigenvector.
\par\parindent=16 pt

Let $A\in{\cal M}_n$ and let $x$ be a nonzero vector of $\IC^n$.
Consider the (unique) representation of $x$ as a sum of generalized
eigenvectors
of $A$,
say,
\begin{center}
$x=x_1+\cdots+x_m$,
\end{center}
\par\parindent=0 pt
where $x_1,\cdots,x_m$ are generalized eigenvector of $A$ corresponding
respectively
to distinct eigenvalues $\lambda_1,\cdots,\lambda_m$.
By definition,
$\rho_{{}_{\scriptstyle x}}(A)$ equals $\max_{1\le i\le
m}\left|\lambda_i\right|$.
By the {\it order of $x$ relative to $A$},
denoted by $\ord_A(x)$,
we mean the maximum of the orders of generalized eigenvectors,
each corresponding to an eigenvalue of modulus $\rho_{{}_{\scriptstyle x}}(A)$,
that appear in the above
representation of $x$.
Now we denote the ordered pair $(\rho_{{}_{\scriptstyle x}}(A), \ord_A(x))$
by $\sp_A(x)$ and refer to it as the {\it spectral pair of $x$ relative to
$A$}.
We also set $\sp_A({\bf 0})=(0, 0)$ to take care of the zero vector ${\bf 0}$.
We denote by $\preceq$ the lexicographic ordering between ordered pairs of real
numbers given by:
\begin{center}
$(\xi_1, \xi_2)\preceq (\eta_1, \eta_2)$ if either $\,\,\xi_1<\eta_1$,
or $\,\,\xi_1=\eta_1$ and $\xi_2\le\eta_2$.
\end{center}
\par\parindent=0 pt

We also write $(\xi_1, \xi_2)\prec (\eta_1, \eta_2)$,
if $(\xi_1, \xi_2)\preceq (\eta_1, \eta_2)$ but the equality does not hold in
the
usual sense.
\par\parindent=16 pt

Using the definition of the spectral pair of a vector (relative to a matrix),
it is easy to establish the following:
\par\parindent=0 pt
\vskip 12 pt

{\it Remark} 4.1.$\,\,$
For any $A\in{\cal M}_n$,
$x, y\in\IC^n$,
$0\not=\lambda\in\IC$,
we have
\par\parindent=16 pt
(i) $\sp_{A}(\lambda x)=\sp_{A}(x)$.
\par\parindent=16 pt
(ii) $\sp_{A}(Ax)$ equals $\sp_A(x)$ if $\rho_{{}_{\scriptstyle x}}(A)>0$,
and equals $(0, \ord_A(x)-1)$ if $\rho_{{}_{\scriptstyle x}}(A)=0$ and $x\not=
{\bf 0}$.
\par\parindent=16 pt
(iii)$\sp_{A}(x+y)\preceq \max\{\sp_{A}(x), \sp_{A}(y)\}$,
where the maximum is taken in the sense of lexicographic ordering.
\par\parindent=16 pt
\vskip 12 pt

Part (iii) of Remark 4.1 is a strengthening of the known result [T--W, Lemma
2.3]
that if $x_1, x_2\in\IC^n$,
then
\begin{center}
$\rho_{{}_{\scriptstyle {x_1+x_2}}}(A)\le\max\{\rho_{{}_{\scriptstyle
{x_1}}}(A), \rho_{{}_{\scriptstyle {x_2}}}(A)\}$.
\end{center}
\par\parindent=0 pt
Indeed,
we have indicated an easy alternative proof.
\par\parindent=16 pt

Here we would like to mention some relevant results which we need later.
For any $n\!\times\! n$ complex matrix $A$ and any eigenvalue $\lambda$ of $A$,
let $E_\lambda^{(0)}$ denote the projection of $\IC^n$ onto the generalized
eigenspace
${\cal N}({(\lambda I-A)}^n)$ along the direct sum of other generalized
eigenspaces of $A$.
Also define the components of $A$ by
\begin{center}
$E_\lambda^{(k)}={(A-\lambda I)}^kE_\lambda^{(0)}, k=0,1,\cdots$.
\end{center}
\par\parindent=0 pt
In [Schn 2, Theorem 5.2 (iii)] it is proved that if $K$ is a proper cone and
$A\in\pi(K)$,
then $\left. E_\rho^{(\nu-1)}\right|_{\,\IR^{\,\,n}}\in\pi(K)$,
where $\rho=\rho(A)$ and $\nu=\nu_{\rho(A)}(A)$.
Furthermore,
we have,
$\rank E_\rho^{(\nu-1)}=\mu_\rho(A)$,
where $\mu_\rho(A)$ is the exponent of $\rho(A)$,
i.e. the number of Jordan blocks of $A$ of maximal order corresponding to
$\rho(A)$.
(This latter fact is always true;
it does not depend on the condition that $A\in\pi(K)$.)
Thus,
$E_\rho^{(\nu-1)}K$ is a cone of dimension $\mu_\rho(A)$ which consists of some
of the distinguished eigenvectors of $A$ corresponding to $\rho(A)$,
together with the zero vector.
For the purpose of later references,
we make a note of the following.
\par\parindent=0 pt
\vskip 12 pt

{\it Remark} 4.2.$\,\,$
If $A\in\pi(K)$,
then there exist $\mu_\rho(A)$ (the exponent of $\rho(A)$)
linearly independent distinguished eigenvectors of $A$ corresponding to
$\rho(A)$,
such that each is associated with a maximal Jordan block,
in the sense that each is a vector of the form ${(A-\rho(A)I)}^{\nu-1}z$
for some generalized eigenvector $z$ of order $\nu=\nu_{\rho(A)}(A)$
corresponding to $\rho(A)$.
\par\parindent=16 pt
\vskip 12 pt

Now we make an important observation.
\par\parindent=0 pt
\vskip 12 pt

{\bf Lemma 4.3.}
{\it Let $A\in\pi(K)$.
For any $x\in\int K$,
we have
\begin{center}
$\sp_{A}(x)=(\rho(A)$,$\nu_{\rho(A)}(A))$.
\end{center}
}
\par\parindent=0 pt
\vskip 12 pt

{\it Proof.}
As we have mentioned above,
$\left. E_\rho^{(\nu-1)}\right|_{\,\IR^{\,\,n}}$ is nonzero and belongs to
$\pi(K)$,
where $\nu=\nu_{\rho(A)}(A)$.
For any $x\in\int K$,
necessarily $E_\rho^{(\nu-1)}x\not= {\bf 0}$;\par\parindent=0 pt
otherwise,
$\left. E_\rho^{(\nu-1)}\right|_{\,\IR^{\,\,n}}$,
and hence $E_\rho^{(\nu-1)}$,
is the zero operator.
Hence,
by definition of $E_\rho^{(\nu-1)}$,
in the representation of $x$ as a sum of generalized eigenvectors of $A$,
there must be a term which is a generalized eigenvector of $A$ corresponding to
$\rho(A)$ and furthermore the order of this generalized eigenvector must be
$\nu$.
Moreover,
by the Perron-Schaefer condition for $A$,
$\nu$ is equal to the maximum of the orders of generalized eigenvectors,
each \par corresponding to an eigenvalue belonging to the peripheral spectrum
of
$A$,
that appear in the representation of $x$.
Therefore,
we have
\begin{center}
\qquad\quad$\sp_{A}(x)=(\rho(A), \nu)$.
\end{center}
\par
\hfill{$\Box$}
\par\parindent=16 pt
\vskip 12 pt

An application of Lemma 4.3 to the restriction of $A$ ($\in\pi(K)$) to the
linear span of an $A$-invariant face $F$ yields immediately the following.
\par\parindent=0 pt
\vskip 12 pt

{\bf Theorem 4.4.}
{\it Let $A\in\pi(K)$,
and let $F$ be an $A$-invariant face of $K$.
Then for any vector $x\in\relint F$,
we have $\sp_A(x)=(\rho_{{}_{\scriptstyle F}},
\nu_{\rho_{{}_{\scriptscriptstyle F}}}(\left. A\right|_{\span F}))$.}
\par\parindent=16 pt
\vskip 12 pt

By Theorem 4.4 and the definition of the spectral pair of a vector,
we have the following.
\par\parindent=0 pt
\vskip 12 pt

{\bf Corollary 4.5.}
{\it Let $A\in\pi(K)$.
If $F$ is an $A$-invariant face of $K$,
then $\rho_{{}_{\scriptstyle x}}(A)=\rho_{{}_{\scriptstyle F}}$ for any
$x\in\relint F$.}
\par\parindent=16 pt
\vskip 12 pt

The result of Theorem 4.4 also says that,
for any $A$-invariant face $F$,
the value of the spectral pair $\sp_A(x)$ is independent of the choice of $x$
from $\relint F$.
Now we are going to show that the same is also true for any face $F$ of $K$.
\par\parindent=0 pt
\vskip 12 pt

{\bf Theorem 4.6.}
{\it Let $A\in\pi(K)$,
and let $F$ be a face of $K$.
For any $x\in\relint F$,
we have
\begin{center}
$\sp_A(x)=(\rho_{{}_{\scriptstyle{\widehat F}}},
\nu_{\rho_{{}_{\scriptscriptstyle{\widehat F}}}}(\left.
A\right|_{\span{\widehat F}}))$,
\end{center}
\par\parindent=0 pt
where $\widehat F$ is the smallest $A$-invariant face of $K$ including $F$.}
\par\parindent=0 pt
\vskip 12 pt

{\it Proof.}
Choose any positive real number $\alpha$ for which $-\alpha$ is not an
eigenvalue of $A$.
Consider any vector $x$ chosen from $\relint F$.
For convenience,
denote the vector ${(\alpha I+A)}^{n-1}x$ by $\widehat x$.
Then the proof of Lemma 2.1 also shows that $\Phi({\widehat x})$ equals
$\widehat F$.
Let
\par\parindent=0 pt
\vskip 10 pt
(1)
\qquad\qquad\qquad\qquad\qquad\quad
$x=x_1+\cdots+x_m$
\vskip 10 pt
\par\parindent=0 pt

where $x_1,\cdots,x_m$ are generalized eigenvectors of $A$ corresponding to
distinct
eigenvalues $\lambda_1,\cdots,\lambda_m$ respectively.
Applying ${(\alpha I+A)}^{n-1}$ to both sides of (1),
we obtain
\par\parindent=0 pt
\vskip 10 pt
(2)
\qquad\qquad\quad\,\,
$\widehat x={(\alpha I+A)}^{n-1}x_1+\cdots+{(\alpha I+A)}^{n-1}x_m$.
\vskip 10 pt
\par\parindent=0 pt
By our choice of $\alpha$,
${(\alpha I+A)}^{n-1}$ is nonsingular;
hence,
each ${(\alpha I+A)}^{n-1}x_i$ is a generalized eigenvector of $A$
corresponding to
$\lambda_i$,
and is of the same order as that of $x_i$.
>From the definition of the spectral pair of a vector,
now it is clear that we have
\begin{center}
$\sp_A(x)=\sp_A(\widehat x)=(\rho_{{}_{\scriptstyle{\widehat
F}}},\nu_{\rho_{{}_{\scriptscriptstyle {\widehat F}}}}(\left.
A\right|_{\span{\widehat F}}))$,
\end{center}
\par\parindent=0 pt
where the last equality follows from Theorem 4.4.
\par
\hfill{$\Box$}
\par\parindent=16 pt
\vskip 12 pt

By Theorem 4.6,
for any face $F$ of $K$,
the value of the spectral pair $\sp_A(x)$ is independent of the choice of $x$
from $\relint F$.
Henceforth,
we shall denote this constant value by $\sp_A(F)$,
and refer to it as the {\it spectral pair of $F$ relative to $A$}.
The spectral pair of a face will prove to be a useful concept,
though it is not an extension of an existing concept for nonnegative matrices.
This probably is due to the fact that its definition implicitly
involves the Perron-Schaefer condition,
a characterizing property of a cone-preserving map.
\par\parindent=16 pt

In passing,
we also take note of the following by-product.
\par\parindent=0 pt
\vskip 12 pt

{\bf Theorem 4.7.}
{\it Let $A\in\pi(K)$.
If ${\bf 0}\not= x\in K$,
then there is always a \par\parindent=0 pt generalized eigenvector $y$ of $A$
corresponding to
$\rho_{{}_{\scriptstyle x}}(A)$ that appears as a term in the representation of
$x$ as a sum of generalized eigenvectors of $A$.\par
Furthermore,
we have $\ord_A(x)=\ord_A(y)$.}
\par\parindent=0 pt
\vskip 12 pt

{\it Proof.}
If the face $\Phi(x)$ is $A$-invariant,
our assertion follows from the proof of Lemma 4.3.
In view of the proof of Theorem 4.6,
our assertion also holds for the case when $\Phi(x)$ is not $A$-invariant.
\par
\hfill{$\Box$}
\par\parindent=16 pt
\vskip 12 pt

The following result should also be clear.
\par\parindent=0 pt
\vskip 12 pt

{\bf Corollary 4.8.}
{\it Let $A\in\pi(K)$,
and let ${\bf 0}\not= x\in K$.
Denote $\ord_A(x)$ by $m$.
Then ${(A-\rho_{{}_{\scriptstyle
x}}(A)I)}^{m-1}E^{(0)}_{\rho_{{}_{\scriptscriptstyle x}}(A)}x$
is a distinguished eigenvector of $A$ \par\parindent=0 pt corresponding to
$\rho_{{}_{\scriptstyle x}}(A)$.}
\par\parindent=16 pt
\vskip 12 pt

In the next result we collect together fundamental properties of spectral pairs
of faces (or,
of vectors in the cone $K$).
\par\parindent=0 pt
\vskip 12 pt

{\bf Theorem 4.9.}
{\it Let $A\in\pi(K)$.}
\par\parindent=16 pt
(i) {\it For any faces $F, G$ of $K$,
we have}
\par\parindent=24 pt
(a) {\it $\sp_A(F)=\sp_A({\widehat F})$,
where $\widehat F$ is the smallest $A$-invariant face of $K$ including $F$.}
\par\parindent=24 pt
(b) {\it If $F\subseteq G$,
then $\sp_A(F)\preceq\sp_A(G)$.}
\par\parindent=24 pt
(c) {\it $\sp_A(F\bigvee G)=\max\{\sp_A(F), \sp_A(G)\}$,
where the maximum is taken in the sense of lexicographic ordering.}
\par\parindent=16 pt
(ii) {\it For any vectors $x, y\in K$,
we have}
\par\parindent=24 pt
(a) {\it $\sp_A(x)=\sp_A({(I+A)}^{n-1}x)$.}
\par\parindent=24 pt
(b) {\it If $x\in\Phi(y)$,
then $\sp_A(x)\preceq\sp_A(y)$.}
\par\parindent=24 pt
(c) {\it $\sp_A(x+y)=\max\{\sp_A(x), \sp_A(y)\}$,
where the maximum is taken in the sense of lexicographic ordering.}
\par\parindent=0 pt
\vskip 12 pt

{\it Proof.}
As can be readily shown,
properties (a),
(b) or (c) of part (i) are each equivalent to the corresponding properties of
part (ii).
So,
for each property,
(a),
(b) or (c),
we need only consider either the one in (i) or the one in (ii).
\par\parindent=16 pt

In view of Theorems 4.4 and 4.6 and the definition of the spectral pair of
a face,
we have $\sp_A(F)=\sp_A({\widehat F})$.
So we need not consider property (a) further.
\par\parindent=16 pt

Let $F, G$ be faces of $K$ such that $F\subseteq G$.
Then clearly we have ${\widehat F}\subseteq {\widehat G}$;
thus,
$\rho_{{}_{\scriptstyle{\widehat F}}}\le\rho_{{}_{\scriptstyle{\widehat G}}}$,
and if $\rho_{{}_{\scriptstyle{\widehat F}}}=\rho_{{}_{\scriptstyle{\widehat
G}}}$
then $\nu_{\rho_{{}_{\scriptscriptstyle{\widehat F}}}}(\left.
A\right|_{\span{\widehat F}})\le
\nu_{\rho_{{}_{\scriptscriptstyle{\widehat G}}}}(\left.
A\right|_{\span{\widehat G}})$.
That is,
\par\parindent=0 pt\hskip 2 pt
we have,
$\sp_A({\widehat F})\preceq\sp_A({\widehat G})$,
and in view of property (a) of (i),
property (b) of (i) follows.
\par\parindent=16 pt

To prove property (c) of (ii),
first note that according to Remark 4.1(iii),
we always have
\begin{center}
$\sp_A(x+y)\preceq\max\{\sp_A(x), \sp_A(y)\}$.
\end{center}
\par\parindent=0 pt
Since $x\in\Phi(x+y)$,
by property (b) of (ii),
$\sp_A(x)\preceq\sp_A(x+y)$.
Similarly,
we also have,
$\sp_A(y)\preceq\sp_A(x+y)$,
and hence
\begin{center}
$\max\{\sp_A(x), \sp_A(y)\}\preceq\sp_A(x+y)$.
\end{center}
\par\parindent=0 pt
Therefore,
the desired equality holds.
\par
\hfill{$\Box$}
\par\parindent=16 pt
\vskip 12 pt

An immediate consequence of property (c) of Theorem 4.9 (i) is that,
if $F$ is an $A$-invariant face that satisfies $\sp_A(F)\succ\sp_A(G)$ for all
$A$-invariant faces $G$ properly included in $F$,
then $F$ is $A$-invariant join-irreducible.
The converse of this result is not true.
For instance,
take $A$ to be the nonnegative matrix
$
\left(
\begin{array}{cc}
1 & 1 \\
0 & 0
\end{array}
\right)
$.
Then $\IR^{\,\,2}_{\,\,+}$ and $\Phi(e_1)$ are the only nonzero $A$-invariant
faces of $\IR^{\,\,2}_{\,\,+}$,
so that $\IR^{\,\,2}_{\,\,+}$ is an $A$-invariant join-irreducible face of
itself.
However,
we have $\sp_A(\IR^{\,\,2}_{\,\,+})=(1,1)=\sp_A(\Phi(e_1))$.
\par\parindent=16 pt
\vskip 12 pt

By Remark 4.1(i),
(ii) and properties (b),
(c) of Theorem 4.9(ii),
we also \par\parindent=0 pt readily obtain the following.
\par\parindent=0 pt
\vskip 12 pt

{\bf Corollary 4.10.}
{\it Let $A\in\pi(K)$.
For any nonnegative real number $\lambda$ and any positive integer $k$,
the set
\begin{center}
$F_{\lambda, k}=\{ x\in K: \sp_A(x)\preceq (\lambda, k)\}$
\end{center}

\par\parindent=0 pt
is an $A$-invariant face of $K$.
Furthermore,
for any nonnegative real numbers $\lambda_1,\lambda_2$,
and positive integers $k_1,k_2$,
if $(\lambda_1,k_1)\preceq (\lambda_2,k_2)$,
then $F_{\lambda_1,k_1}\subseteq F_{\lambda_2,k_2}$.}
\par\parindent=16 pt
\vskip 12 pt

For a similar reason,
for any positive real number $\lambda$,
the set
\begin{center}
$F_\lambda=\{ x\in K: \rho_{{}_{\scriptstyle x}}(A)<\lambda\}$
\end{center}
\par\parindent=0 pt
is an $A$-invariant face of $K$;
the greater $\lambda$ is,
the larger is the face $F_\lambda$.
In the case of a positive linear operator $T$ defined on a real Banach lattice
$E$,
Meyer-Nieberg (see [MN 1] or [MN 2, pp. $\!\!$293])
has discovered a class of $T$-invariant ideals $J_\lambda$,
which is precisely the linear span of our $F_\lambda$ when the underlying
Banach
lattice is $\IR^{\,\,n}$ (with positive cone $\IR^{\,\,n}_{\,\,+}$).
We are going to explain.
\par\parindent=16 pt

For each real number $\lambda>0$,
the ideal $J_\lambda$ introduced by Meyer-Nieberg is given by:
\begin{center}
{$J_\lambda=\{ x\in E:$ the series $\sum^\infty_{j=0}T^j\left| x\right| /
\lambda^{j+1}$
is convergent$\}$.}
\end{center}
\par\parindent=0 pt
Here we follow the standard notation of Banach lattice theory (for reference,
see [Scha]) and use $\left| x\right|$ to denote the modulus of $x$.
First,
we show that
\begin{center}
$F_\lambda=\{ x\in\IR^{\,\,n}_{\,\,+}:$ the series $\sum^\infty_{j=0}T^jx/
\lambda^{j+1}$
is convergent$\}$.
\end{center}
\par\parindent=0 pt
If $x\in F_\lambda$,
then by definition,
$\rho_{{}_{\scriptstyle x}}(T)<\lambda$.
But $\rho_{{}_{\scriptstyle x}}(T)=\rho(\left. T\right|_{W_x})$
(see [T--W, Theorem 2.2]),
where $W_x$ is the cyclic space generated by $x$,
so ${(\lambda I-\left. T\right|_{W_x})}^{-1}$ exists and is given by:
\begin{center}
${(\lambda I-\left. T\right|_{W_x})}^{-1}=\lambda^{-1}
{{{\scriptstyle\infty}\atop{\textstyle\sum}}\atop{j=0}} {(\left.
T\right|_{W_x}/\lambda)}^j$.
\end{center}
\par\parindent=0 pt
Hence,
the series $\sum^\infty_{j=0}T^jx/\lambda^{j+1}=\sum^\infty_{j=0}{(\left.
T\right|_{W_x})}^jx/\lambda^{j+1}$
is convergent.
\par\parindent=0 pt Conversely,
if the series $\sum^\infty_{j=0}T^jx/\lambda^{j+1}$ converges,
then
\begin{center}
${{}\atop{{\textstyle\lim}\atop{\scriptstyle{m\to\infty}}}}{(T/\lambda)}^m(T^ix
)={\bf 0}\,\,$ for $\,\, i=0,1,2,\cdots$,
\end{center}
\par\parindent=0 pt
and hence
\begin{center}
${{}\atop{{\textstyle\lim}\atop{\scriptstyle{m\to\infty}}}}{(T/\lambda)}^my={\b
f 0}\,\,$ for every $\,\, y\in W_x$.
\end{center}
\par\parindent=0 pt
It follows that $\lim_{m\to\infty}{(\left. T\right|_{W_x}/\lambda)}^m={\bf 0}$,
and so $\rho_{{}_{\scriptstyle x}}(T)<\lambda$.
This establishes our claim on the set $F_\lambda$.
\par\parindent=16 pt

For any $x\in\IR^{\,\,n}$,
denote the positive part and the negative part of $x$ \par\parindent=0 pt
respectively by $x^+$ and $x^-$.
Let $x\in\span F_\lambda$.
Since $F_\lambda$ is a face of $\IR^{\,\,n}_{\,\,+}$,
it is clear that $x^+$ and $x^-$ both belong to $F_\lambda$.
Hence,
both of the series $\sum^\infty_{j=0}T^jx^+/\lambda^{j+1}$ and
$\sum^\infty_{j=0}T^jx^-/\lambda^{j+1}$
converge,
and so does the series \par $\sum^\infty_{j=0}T^j\left| x\right|
/\lambda^{j+1}$,
as $\left| x\right| =x^++x^-$.
This proves the inclusion $\span F_\lambda\subseteq J_\lambda$.
Conversely,
if $x\in J_\lambda$,
then it is clear that both of the series
$\sum^\infty_{j=0}T^jx^+/\lambda^{j+1}$
and $\sum^\infty_{j=0}T^jx^-/\lambda^{j+1}$ converge;
hence,
$x^+, x^-\!\in F_\lambda$,
and so $x=x^+-x^-\in\span F_\lambda$.
This proves the reversed inclusion and hence the desired equality.
\par\parindent=16 pt

So,
in a sense,
our Corollary 4.10 is an extension and also a refinement of the above-mentioned
result of Meyer-Nieberg for the finite-dimensional case.
\par\parindent=16 pt
\vskip 9 pt

Each property about spectral pairs given in Theorem 4.9 clearly implies a
corresponding
property about (local) spectral radii.
For instance,
by property (c) of (i) we have the following:
{\it If $F, G$ are $A$-invariant faces of $K$,
then $\rho_{{}_{\scriptstyle {F\bigvee G}}}=\max\{\rho_{{}_{\scriptstyle F}},
\rho_{{}_{\scriptstyle G}}\}$.}
We leave it to the reader to formulate the details.
\par\parindent=16 pt
\vskip 2 pt

Now we give equivalent conditions for an $A$-invariant face of $K$ to be
distinguished.
\par\parindent=0 pt
\vskip 12 pt

{\bf Theorem 4.11.}
{\it Let $A\in\pi(K)$,
and let $F$ be a face of $K$.
The following conditions are equivalent:}
\par\parindent=16 pt
(a) {\it $F$ is a distinguished $A$-invariant face of $K$.}
\par\parindent=16 pt
(b) {\it $F$ is an $A$-invariant face,
there is a unique (up to positive multiples) eigenvector of $A$ in $F$ that
corresponds
to $\rho_{{}_{\scriptstyle F}}$,
and this eigenvector lies in $\relint F$.}
\par\parindent=16 pt
(c) {\it $F=\Phi(x)$ for some extremal distinguished eigenvector of $A$.}
\par\parindent=16 pt
(d) {\it $F$ is an $A$-invariant join-irreducible face of $K$,
which contains in its relative interior an eigenvector of $A$.}
\par\parindent=0 pt
\vskip 12 pt

{\it Proof.}
\par\parindent=16 pt
(a) $\Longrightarrow$ (b):
Let $x$ be an eigenvector of $A$ corresponding to $\rho_{{}_{\scriptstyle F}}$
that lies in $F$.
If $x\notin\relint F$,
then $\Phi(x)$ is a nonzero $A$-invariant face properly included in
$F$ such that $\rho_{{}_{\scriptstyle {\Phi(x)}}}=\rho_{{}_{\scriptstyle F}}$.
This contradicts the assumption that $F$ is a distinguished $A$-invariant face.
If $x\in\relint F$ and there is an eigenvector,
not a multiple of $x$,
corresponding to $\rho_{{}_{\scriptstyle F}}$ that lies in $F$,
then we readily obtain an eigenvector corresponding to $\rho_{{}_{\scriptstyle
F}}$
that lies in $\rbd F$,
and by the first part of our argument we also obtain a contradiction.
\par\parindent=16 pt
The equivalence of (b) and (c) is clear.
\par\parindent=16 pt
(b) $\Longrightarrow$ (a):
Suppose that there is a nonzero $A$-invariant face $G$ properly included in $F$
such that $\rho_{{}_{\scriptstyle G}}=\rho_{{}_{\scriptstyle F}}$.
Then there exists an eigenvector of $A$ corresponding to
$\rho_{{}_{\scriptstyle F}}$
that lies in $\,\rbd F$.
This contradicts (b).
\par\parindent=16 pt
(a) $\Longrightarrow$ (d):
By condition (b),
which is already proved to be equivalent to condition (a),
$A$ has an eigenvector that lies in $\relint F$.
Let $F_1$,
$F_2$ be $A$-invariant faces such that $F=F_1\bigvee F_2$.
Then $\rho_{{}_{\scriptstyle F}}=\max\{\rho_{{}_{\scriptstyle {F_1}}},
\rho_{{}_{\scriptstyle {F_2}}}\}=\rho_{{}_{\scriptstyle {F_1}}}$,
say.
Since $F$ is distinguished $A$-invariant,
it follows that $F_1=F$.
This proves the $A$-invariant join-irreducibility of $F$.
\par\parindent=16 pt
(d) $\Longrightarrow$ (c):
By condition (d) $A$ has an eigenvector $x$ that lies in $\relint F$
(and corresponds to $\rho_{{}_{\scriptstyle F}}$).
Assume to the contrary that $x$ is not an extremal distinguished eigenvector.
Then there exist eigenvectors $x_1$,
$x_2$ of $A$ \par\parindent=0 pt corresponding to $\rho_{{}_{\scriptstyle F}}$
such that both vectors
belong to $\rbd F$ and $x_1+x_2=x$.
Then $\Phi(x_1)$ and $\Phi(x_2)$ are both $A$-invariant faces of $K$ properly
included in $F$ such that $\Phi(x_1)\bigvee\Phi(x_2)=\Phi(x_1+x_2)=F$.
This contradicts the assumption that $F$ is $A$-invariant join-irreducible.
\par
\hfill{$\Box$}
\par\parindent=16 pt
\vskip 12 pt

By Theorem 3.6(iv),
if $P$ is a nonnegative matrix and $I$ is an initial subset for $P$,
then $I$ is determined by a semi-distinguished (respectively,
distinguished) class if and only if $F_I$ is a $P$-invariant join-irreducible
face which contains in its relative interior a generalized eigenvector
(respectively,
an eigenvector) of $P$.
And in view of the equivalent condition (d) of Theorem 4.11 for a distinguished
$A$-invariant face,
given any $A\in\pi(K)$,
we shall call a face $F$ of $K$ a {\it semi-distinguished $A$-invariant face}
(associated with $\lambda$) if $F$ is an $A$-invariant join-irreducible face
which
contains in its relative interior a generalized eigenvector of $A$
(corresponding to $\lambda$).
Clearly a distinguished $A$-invariant face is semi-distinguished A-invariant.
\par\parindent=16 pt

Note that we define a nonzero $A$-invariant face $F$ to be distinguished by the
property that $\rho_{{}_{\scriptstyle G}}<\rho_{{}_{\scriptstyle F}}$ for any
nonzero $A$-invariant face $G$ properly \par\parindent=0 pt included in $F$.
Here a semi-distinguished $A$-invariant face is not defined to be one that
satisfies the weaker property obtained by replacing the strict inequality by
the
weak inequality of the same type,
because this latter \par property is always satisfied!
Another reasonable defining property for a semi-distinguished $A$-invariant
face is the following:
\par\parindent=16 pt
\vskip 12 pt

{\it $F$ is nonzero $A$-invariant,
and $\sp_A(G)\prec\sp_A(F)$ for any $A$-invariant face $G$ properly included in
$F$.}
\par\parindent=0 pt
\vskip 12 pt

In view of Theorem 3.6,
it is easy to see that in the case when $A$ is a \par\parindent=0 pt
nonnegative matrix,
the preceding property is also another equivalent \par condition for $F$ to be
a
semi-distinguished $A$-invariant face.
Indeed,
as we shall see,
the same is also true whenever the underlying cone $K$ has the property
that the dual cone of each of its faces is facially exposed
(which is the case if $K$ is polyhedral),
but is not so for a general proper cone $K$.
\par\parindent=16 pt
\vskip 12 pt

We shall need the following general result,
which is a consequence of the existence of a Jordan basis for a matrix.
\par\parindent=0 pt
\vskip 12 pt

{\it Remark} 4.12.$\,\,$
Let $A\in {\cal M}_n$ be such that $\rho(A)$ is an eigenvalue.
If $z$ is an eigenvector of $A^T$ corresponding to $\rho(A)$ which is
associated with a maximal Jordan block,
then $z$ is orthogonal to every vector $x$ for which
$\sp_A(x)\prec (\rho(A),\nu_{\rho(A)}(A))$.
\par\parindent=0 pt
\vskip 25 pt

\begin{center}
5. A{\scriptsize N EXTENSION

OF THE ROTHBLUM INDEX THEOREM}
\end{center}
\par\parindent=16 pt

Hershkowitz and Schneider (see [H--S, the paragraph following Theorem 5.10])
proved that the famous Rothblum Index theorem (i.e. Theorem $A$ of Section 1)
can be extended to cover the case of a distinguished eigenvalue $\lambda$,
provided that $\lq\lq$ basic classes $"$ are replaced by
$\lq\lq$ semi-distinguished classes associated with $\lambda$ $"$
and $\nu$ is replaced by $m_\lambda$,
the maximal order of \par\parindent=0 pt distinguished generalized eigenvectors
of $P$ corresponding to $\lambda$. This is indeed an extension,
because for a nonnegative matrix $P$,
as an outcome of the Nonnegative Basis theorem,
we have $m_{\rho(P)}=\nu_{\rho(P)}(P)$.
\par\parindent=16 pt

In the case of a linear mapping preserving a polyhedral cone,
we have the following extension of the Rothblum Index theorem.
\par\parindent=0 pt
\vskip 12 pt

{\bf Theorem 5.1.}
{\it Let $K$ be a polyhedral cone,
and let $A\in\pi(K)$.
Let $\lambda$ be a distinguished eigenvalue of $A$ for $K$.
Denote by $m_\lambda$ the maximal order of distinguished generalized
eigenvectors
of $A$ corresponding to $\lambda$.
Then there is a chain of $m_\lambda$ distinct semi-distinguished $A$-invariant
faces of $K$ associated with $\lambda$,
but there is no such chain with more than $m_\lambda$ members.}
\par\parindent=16 pt
\vskip 12 pt

As can be readily seen,
the proof of Theorem 5.1 should consist of two parts:
to construct the desired chain of semi-distinguished $A$-invariant faces,
and to establish the maximality of its cardinality.
The proof of the \par\parindent=0 pt maximality of cardinality depends on the
following result which has interest
of its own.
\par\parindent=0 pt
\vskip 12 pt

{\bf Lemma 5.2.}
{\it Let $K$ be a polyhedral cone,
and let $A\in\pi(K)$.
If $K$ is a semi-distinguished $A$-invariant face of itself,
then $\sp_A(K)\succ\sp_A(F)$ for any $A$-invariant face $F$ properly included
in $K$.}
\par\parindent=0 pt
\vskip 12 pt

{\it Proof.}
Since $\int K$ contains a generalized eigenvector of $A$,
by [Tam 3, \par\parindent=0 pt Theorem 5.1]
$\rho(A)$ is the only distinguished eigenvalue of $A^T$ for $K^*$.
\par\parindent=16 pt

We contend that the cone ${\cal N}(\rho(A)I-A^T)\bigcap K^*$ is of dimension
one.
Suppose not.
Choose distinct extreme vectors $z_1, z_2$ from this cone.
If $\Phi(z_1)\bigwedge\Phi(z_2)\not=\{ {\bf 0}\}$,
then being a nonzero $A^T$-invariant face $\Phi(z_1)\bigwedge\Phi(z_2)$ would
contain an eigenvector of $A^T$,
which necessarily
corresponds to $\rho(A)$,
and hence is a multiple of $z_1$,
and also of $z_2$,
which is a contradiction.
Hence,
we have
$\Phi(z_1)\bigwedge\Phi(z_2)=\{ {\bf 0}\}$.
By our hypothesis,
$K^*$ is polyhedral.
So each face of $K^*$ is exposed and we can write $\Phi(z_i)$
($i=1,2$) as $d_K\circ d_{K^*}(\Phi(z_i))$.
Thus we have
\begin{eqnarray*}
0 & = & \Phi(z_1)\bigwedge\Phi(z_2)					 \\
  & = & d_K\circ d_{K^*}(\Phi(z_1))\bigwedge d_K\circ d_{K^*}(\Phi(z_2)) \\
  & = & d_K\left(d_{K^*}(\Phi(z_1))\bigvee d_{K^*}(\Phi(z_2))\right),
\end{eqnarray*}
which implies that $K$ is the join of the nontrivial $A$-invariant faces
$d_{K^*}(\Phi(z_1))$ and $d_{K^*}(\Phi(z_2))$.
So again we arrive at a contradiction.
This establishes our contention.
And since $\rho(A)$ is the only distinguished eigenvalue of $A^T$,
$A^T$ has (up to multiples) a unique distinguished eigenvector.
\par\parindent=16 pt

In view of Remark 4.2,
another consequence of our contention is that,
$A^T$ has precisely one maximal Jordan block corresponding to $\rho(A)$.
Furthermore,
the unique distinguished eigenvector,
say $w$,
of $A^T$ for $K^*$ must correspond to this maximal block.
Thus,
we have
\begin{center}
${(\span\{ w\})}^\bot=\bigoplus_{\lambda\in\sigma(A)\backslash\{\rho(A)\}
}{\cal N}({(\lambda I-A)}^n)\bigoplus U$,
\end{center}
\par\parindent=0 pt
where $U$ is the space of generalized eigenvectors of $A$ corresponding to
$\rho(A)$
of order less than or equal to $\nu_{\rho(A)}(A)-1$,
together with the zero vector.
\par\parindent=16 pt

Now let $F$ be a proper nonzero $A$-invariant face of $K$.
Since the $A^T$-invariant face $d_K(F)$ of $K^*$ must contain $w$,
we have $\span F\subseteq {(\span\{ w\})}^\bot$.
In view of the above direct decomposition for ${(\span\{ w\})}^\bot$,
it follows that if $\span F$ contains a generalized eigenvector of $A$
corresponding to $\rho(A)$,
then the order of such generalized eigenvector must be less than or equal to
$\nu_{\rho(A)}(A)-1$. So we have
\begin{center}
\qquad\qquad $\sp_A(F)\preceq (\rho(A), \nu_{\rho(A)}(A)-1)\prec\sp_A(K)$.
\par
\hfill{$\Box$}
\end{center}
\par\parindent=0 pt
\vskip 12 pt

{\it Proof of Theorem 5.1.}
By definition of $m_\lambda$ there exists a distinguished \par\parindent=0 pt
generalized
eigenvector $x$ of $A$ corresponding to $\lambda$ of order $m_\lambda$.
Clearly the vector $y={(I+A)}^{n-1}x$ is also a distinguished generalized
eigenvector \par corresponding to $\lambda$ of order $m_\lambda$,
and by Lemma 2.1 $\Phi(y)$ is an $A$-invariant face.
Thus,
it is possible to find an $A$-invariant face of $K$ that contains in its
relative
interior a generalized eigenvector of $A$ corresponding to $\lambda$ of order
$m_\lambda$.
Choose $F_{m_\lambda}$ to be an $A$-invariant face of $K$ minimal with respect
to this property.
By definition,
clearly we have $\sp_A(F_{m_\lambda})=(\lambda, m_\lambda)$.
We contend \par that $F_{m_\lambda}$ is $A$-invariant join-irreducible,
and hence is semi-distinguished.
Suppose not.
Then we can find $A$-invariant faces $G_1$,
$G_2$ properly included in $F_{m_\lambda}$ such that $F_{m_\lambda}=G_1\bigvee
G_2$.
By property (c) of Theorem 4.9(i),
either \par $\sp_A(G_1)$ or $\sp_A(G_2)$ is equal to $\sp_A(F_{m_\lambda})$,
say,
$\sp_A(G_1)$.
Then by definition of the spectral pair of a face and by Theorem 4.4,
$\rho_{{}_{\scriptstyle {G_1}}}=\lambda$ and $\nu_\lambda(\left.
A\right|_{\span G_1})=m_\lambda$.
Since $G_1$ is polyhedral,
according to [Tam 3, Theorem 7.5(ii)],
we can find in $G_1$ a generalized eigenvector $u$ of $A$ corresponding to
$\lambda$ of order $m_\lambda$ such that ${(A-\lambda I)}^iu\in G_1$ for
$i=1,\cdots,m_\lambda-1$.
Thus,
$\Phi({(I+A)}^{n-1}u)$ is an $A$-invariant face included in $G_1$,
and hence properly included in $F_{m_\lambda}$,
with the property that it contains in its relative interior a generalized
eigenvector
of $A$ corresponding to $\lambda$ of order $m_\lambda$.
This contradicts the minimality of $F_{m_\lambda}$,
and hence proves our contention for $F_{m_\lambda}$.
\par\parindent=16 pt

Next,
in view of the polyhedrality of the cone $F_{m_\lambda}$,
by [Tam 3, Theorem 7.5(ii)] again,
we can find in $F_{m_\lambda}$ a generalized eigenvector of $A$ corresponding
\par\parindent=0 pt to $\lambda$
of order $m_\lambda-1$.
Then we can repeat the above argument to obtain a semi-distinguished
$A$-invariant
face $F_{m_\lambda-1}$ of $F_{m_\lambda}$ such that
$\sp_A(F_{m_\lambda-1})=(\lambda, m_\lambda-1)$.
Continuing in this way,
after a finite number of steps,
we \par eventually arrive at our desired chain of $A$-invariant faces.
\par\parindent=16 pt

To show that our chain of $A$-invariant faces is of maximum cardinality,
let $H_1\subseteq\cdots\subseteq H_p$ be an arbitrary chain of distinct
semi-distinguished $A$-invariant faces associated with $\lambda$.
Then for each $i$,
$1\le i\le p$,
$\rho_{{}_{\scriptstyle {H_i}}}=\lambda$.
In view of Lemma 5.2 and the definition of $m_\lambda$,
it follows that we have
\begin{center}
$(\lambda, 1)\preceq\sp_A(H_1)\prec\sp_A(H_2)\prec\cdots\prec\sp_A(H_p)\preceq
(\lambda, m_\lambda)$.
\end{center}
\par\parindent=0 pt
Hence,
$p\le m_\lambda$.
This shows that the chain of $A$-invariant faces (with the desired property)
we have constructed is of maximum cardinality.
\par
\hfill{$\Box$}
\par\parindent=16 pt
\vskip 12 pt

When $K$ is a nonpolyhedral cone,
the maximum cardinality of a chain of semi-distinguished $A$-invariant faces of
$K$
associated with a distinguished eigenvalue $\lambda$ can be less than,
equal to,
or greater than $m_\lambda$,
where $m_\lambda$ has the same meaning as before.
The following examples will serve to illustrate this point and also other
related
points.
\par\parindent=0 pt
\vskip 12 pt

{\bf Example 5.3.}$\,\,$
Let $\alpha$ be a given real number with $0<\alpha <1$.
Let $C$ be the closed convex set in $\IR^{\,\, 2}$ with extreme points
${(k, {\alpha}^{k-1})}^T$,
$k=1,2,\cdots$,
and with recession cone ${\bf 0}^+C=\{\lambda{(1, 0)}^T: \lambda\ge 0\}$.
Let $K$ be the proper cone of $\IR^{\,\, 3}$ given by:
\begin{center}
$K=\{\lambda {x\choose 1}\in\IR^{\,\, 3}: x\in C, \lambda\ge 0\}\,
\bigcup\,\{\lambda {(1,0,0)}^T: \lambda\ge 0\}$.
\end{center}
\par\parindent=0 pt
Let
\begin{center}
$A=
\left[
\begin{array}{ccc}
1 & 0	   & 1 \\
0 & \alpha & 0 \\
0 & 0	   & 1
\end{array}
\right]
$
.
\end{center}
\par\parindent=0 pt
(The cone $K$ and the matrix $A$ have been considered in [T--W, Example 3.7].)
\par\parindent=0 pt
Straightforward calculations show that we have
$A{(k, \alpha^{k-1}, 1)}^T={(k+1, \alpha^k, 1)}^T$ for $k=1,2,\cdots$,
and $A{(1,0,0)}^T={(1,0,0)}^T$.
It follows that $A\in\pi(K)$.
Clearly $\rho(A)=1$ and $\nu_{\rho(A)}(A)=2$.
It is also not difficult to see that $A$ has only two nonzero $A$-invariant
faces,
namely,
$K$ and $F$,
where $F=\Phi({(1,0,0)}^T)$,
both of which are $A$-invariant join-irreducible,
but only $\Phi({(1,0,0)}^T)$ is semi-distinguished $A$-invariant.
In this case,
we have $m_{\rho(A)}=1<2=\nu_{\rho(A)}(A)$.
Here are two other observations we can make of this example
(for the nonpolyhedral case):
\par\parindent=16 pt
(i)
There does not exist a semi-distinguished $A$-invariant face $G$ of $K$
such that $\sp_A(G)=\sp_A(K)$.
Indeed,
$F$ is the only semi-distinguished $A$-invariant face of $K$,
and we have $\sp_A(F)=(1,1)\prec (1,2)=\sp_A(K)$.
\par\parindent=16 pt
(ii)
$K$ is not a semi-distinguished $A$-invariant face of itself,
though it has the property that $\sp_A(K)\succ\sp_A(G)$ for all $A$-invariant
faces $G$ properly included in $K$.

\par\parindent=0 pt
\vskip 12 pt

{\bf Example 5.4.}$\,\,$
Let $C$ be the closed convex set in $\IR^{\,\, 2}$ with recession cone
${\bf 0}^+C=\{\lambda{(1,0)}^T: \lambda\ge 0\}$,
and with extreme points ${(k(k-1)/2, k)}^T$,
$k=0,1,2,\cdots$.
Let $K$ be the proper cone of $\IR^{\,\, 3}$  given by:
\begin{center}
$K=\{\lambda {x\choose 1}\in\IR^{\,\, 3}: x\in C, \lambda\ge 0\}\,
\bigcup\,\{ {x\choose 0}\in\IR^{\,\, 3}: x\in {\bf 0}^+C\}$.
\end{center}
\par\parindent=0 pt
Let $A=J_3(1)$,
the $3\!\times\! 3$ (upper triangular) elementary Jordan matrix
\par\parindent=0 pt corresponding
to the eigenvalue $1$.
(The matrix $A$ and the cone $K$ have been considered in [Tam 3, Example 5.5]).
It is straightforward to check that $A{(1,0,0)}^T={(1,0,0)}^T$,
and $A({(k(k-1)/2, k, 1)}^T)={(k(k+1)/2, k+1, 1)}^T$ for $k=0,1,2,\cdots$.
It follows that $A\in\pi(K)$.
Clearly,
$A$ has only two nonzero $A$-invariant faces,
namely $K$ and $\Phi({(1,0,0)}^T)$,
both of which are semi-distinguished $A$-invariant.
In this case,
every nonzero vector of $K$ is a distinguished generalized eigenvector of $A$
corresponding to $\rho(A)$ ($=1$).
Also $m_{\rho(A)}=\nu_{\rho(A)}(A)=3$,
but it is impossible to find a chain of three distinct semi-distinguished
$A$-invariant faces of $K$.
\par\parindent=0 pt
\vskip 12 pt

{\bf Example 5.5.}$\,\,$
Let $K$ be a $3$-dimensional proper cone whose intersection with the hyperplane
$\{ {(\xi_1, \xi_2, \xi_3)}^T: \xi_1+\xi_2+\xi_3=1\}$
is a convex set $C$ lying inside the simplex $\conv\{ e_1, e_2, e_3\}$,
where $e_1$,
$e_2$,
$e_3$ are the standard unit vectors of $\IR^{\,\, 3}$,
and with the relative boundary formed by the line segment joining $e_1$ and
$e_2$,
together with some strictly convex smooth curve.
\par\parindent=16 pt

Let $A=J_2(0)\bigoplus J_1(0)$.
(Note that,
in fact,
$A=e_1{e_2}^T$.)
It is straight-forward to check that $A\in\pi(K)$ and also that $K$ has three
nonzero $A$-invariant faces,
which are pairwise comparable,
namely,
$\Phi(e_1)$,
$\Phi(e_1+e_2)$ and $K$ itself.
Indeed,
each of them is a semi-distinguished $A$-invariant face associated with
$\rho(A)$
($=0$).
So,
in this case,
the cardinality of a maximum chain of semi-distinguished $A$-invariant
faces associated with $\rho(A)$ is $3$ and is greater than $m_{\rho(A)}$,
which is $2$.
Note that we also have
\begin{center}
$\sp_A(K)=(0,2)=\sp_A(\Phi(e_1+e_2))$,
\end{center}
\par\parindent=0 pt
even though the $A$-invariant face $\Phi(e_1+e_2)$ is properly included in
the semi-distinguished $A$-invariant face $K$.
\par\parindent=16 pt
It is also worth noting that in this example the fact that $\rho(A)=0$
is not a crucial condition.
If we replace $A$ by $A+I$,
the example still works.
}

{\baselineskip	0.6 cm
\par\parindent=0 pt
\vskip 25 pt

\begin{center}
6. S{\scriptsize EMI-DISTINGUISHED $A$-INVARIANT FACES}
\end{center}
\par\parindent=16 pt

Certainly,
a natural question to ask is,
when a nonzero $A$-invariant face $F$ of $K$ is semi-distinguished
$A$-invariant.
Clearly $F$ is a semi-distinguished $A$-invariant face of $K$ if and only if
$F$
is a semi-distinguished $\left. A\right|_{\span F}$-invariant face of itself.
Thus,
it is fundamental to consider the problem of when $K$ is a semi-distinguished
$A$-invariant face of
itself.
We have the following result which is motivated by the proof of Lemma 5.2.
\par\parindent=0 pt
\vskip 12 pt

{\bf Theorem 6.1.}
{\it Let $K$ be a proper cone and let $A\in\pi(K)$.
Consider the following conditions:}
\par\parindent=16 pt
(a)
{\it $A^T$ has
(up to multiples)
a unique distinguished eigenvector for $K^*$.}
\par\parindent=16 pt
(b)
{\it $K$ is a semi-distinguished $A$-invariant face of itself.}
\par\parindent=16 pt
(c)
{\it $\sp_A(K)\succ\sp_A(F)$ for all $A$-invariant faces $F$ properly included
in $K$.}
\par\parindent=0 pt
{\it We always have} (a)$\Longrightarrow$(b) {\it and} (a)$\Longrightarrow$(c).
{\it If $K^*$ is a facially exposed cone,
then conditions} (a),
(b) {\it and} (c) {\it are equivalent.}
\par\parindent=0 pt
\vskip 12 pt

{\it Proof.}
\par\parindent=16 pt
(a)$\Longrightarrow$(b):
When $K^*$ contains a unique eigenvector of $A^T$,
clearly $\rho(A)$ is the only distinguished eigenvalue of $A^T$ for $K^*$.
In this case,
by [Tam 3, Theorem 5.1] $K$ must contain in its interior a generalized
eigenvector of $A$.
So it suffices to establish the $A$-invariant join-irreducibility of $K$.
Assume to the contrary that there exist $A$-invariant faces $F_1, F_2$ of $K$,
both different from $K$,
such that $K=F_1\bigvee F_2$.
Then $d_K(F_1)$,
$d_K(F_2)$ are nonzero $A^T$-invariant faces of $K^*$,
and hence must contain distinguished eigenvectors of $A^T$,
say,
$z_1$ and $z_2$ respectively.
But
\begin{center}
$d_K(F_1)\bigwedge d_K(F_2)=d_K(F_1\bigvee F_2)=d_K(K)=\{ {\bf 0}\}$,
\end{center}
\par\parindent=0 pt
so the vectors $z_1$ and $z_2$ are linearly independent.
This contradicts condition (a).
\par\parindent=16 pt

(a)$\Longrightarrow$(c):
Follows from the argument used in the second half of the proof of Lemma 5.2.
\par\parindent=16 pt

Now suppose,
in addition,
that $K^*$ is a facially exposed cone.
\par\parindent=16 pt

(b)$\Longrightarrow$(a):
Follows from the argument used in the first half of the proof of Lemma 5.2.
\par\parindent=16 pt

(c)$\Longrightarrow$(b):
By condition (c)
$K$ is an $A$-invariant join-irreducible face of itself.
Suppose that $K$ is not semi-distinguished $A$-invariant.
Then \par\parindent=0 pt necessarily $\int K$ contains no generalized
eigenvector of $A$.
Hence,
$A^T$ has a distinguished eigenvalue for $K^*$ other than $\rho(A)$,
say $\lambda$,
and let $z\in K^*$ be a corresponding eigenvector.
Choose a distinguished eigenvector of $A^T$ corresponding to $\rho(A)$
that is associated with a maximal block,
say $w$.
We contend that $w\in d_K\circ d_{K^*}(\Phi(z))$.
Take any vector $x$ from $d_{K^*}(\Phi(z))$.
Since $z$ is an eigenvector of $A^T$,
$d_{K^*}(\Phi(z))$ is an $A$-invariant face of $K$,
different from $K$.
In view of condition (c),
$\sp_A(x)\preceq\sp_A(d_{K^*}(\Phi(z)))\prec\sp_A(K)=(\rho(A),\nu_{\rho(A)}(A))
$;
so by Remark 4.12
we have $\langle x,w\rangle=0$.
This proves our contention.
Note that $w\not\in\Phi(z)$,
as $\sp_{A^T}(w)=(\rho(A),1)\not\prec (\lambda,1)=\sp_{A^T}(z)$.
Hence,
we have $\Phi(z)\not= d_K\circ d_{K^*}(\Phi(z))$.
In other words,
$\Phi(z)$ is not an exposed face.
This contradicts our hypothesis that $K^*$ is a facially exposed cone.
\par
\hfill{$\Box$}
\par\parindent=16 pt
\vskip 12 pt

Our proof of Theorem 6.1,
(c)$\Longrightarrow$(b) actually also establishes the following positive
results:
\par\parindent=0 pt
\vskip 12 pt

{\bf Theorem 6.2.}
{\it Let $K$ be a proper cone and let $A\in\pi(K)$.
Suppose that $\sp_A(K)\succ\sp_A(F)$ for all $A$-invariant faces $F$ properly
included in $K$.
If $K$ is not a semi-distinguished $A$-invariant face of itself,
then for any distinguished eigenvector $z$ of $A^T$ for $K^*$ corresponding
to a distinguished eigenvalue
other than $\rho(A)$,
$\Phi(z)$ is a non-exposed face of $K^*$.}
\par\parindent=0 pt
\vskip 12 pt

{\bf Corollary 6.3.}
{\it Let $A\in\pi(K)$,
and suppose that $\sp_A(K)\succ\sp_A(F)$ for all $A$-invariant faces $F$
properly included in $K$.
If we have,
for any non-exposed face $G$ of $K^*$,
the face $d_{K^*}(G)$ of $K$ is not $A$-invariant,
then $K$ is a semi-distinguished $A$-invariant face of itself.}
\par\parindent=16 pt
\vskip 12 pt

It is worth noting that condition (a) of Theorem 6.1 is slightly
weaker than the following equivalent condition for $A^T$ to be
$K^*$-irreducible:
\par\parindent=16 pt

{\it $A^T$ has
(up to multiples)
a unique distinguished eigenvector for $K^*$,
and this eigenvector lies in $\int K^*$.}
\par\parindent=16 pt

In general,
conditions (b),
(c) of Theorem 6.1 are logically independent.
By Examples 5.3 and 5.5 respectively,
we do not have the implications (c)$\Longrightarrow$(b)
and (b)$\Longrightarrow$(c)
(and hence also not (c)$\Longrightarrow$(a) and (b)$\Longrightarrow$(a)).
\par\parindent=16 pt

The following result clarifies the logical relations between a number of
natural conditions
on a cone-preserving map that involve the concepts of
semi-distinguished $A$-invariant faces,
distinguished generalized eigenvectors,
and spectral pairs of faces.
\par\parindent=0 pt
\vskip 12 pt

{\bf Theorem 6.4.}
{\it Let $A\in\pi(K)$,
and let $\lambda$ be a distinguished eigenvalue of $A$ for $K$.
Consider the following conditions:}
\par\parindent=16 pt
(a)
\par\parindent=24 pt
(i)
{\it Any $A$-invariant face of $K$ which contains in its relative interior a
generalized eigenvector of $A$
corresponding to $\lambda$ can be expressed as a join of semi-distinguished
$A$-invariant
faces associated with $\lambda$.}
\par\parindent=24 pt
(ii)
{\it Let $C$ denote the cone $K\bigcap {\cal N}({(\lambda I-A)}^n)$,
and let $B$ be the \par\parindent=0 pt restriction map $\left. A\right|_{\span
C}$.
For any semi-distinguished $B$-invariant face $F$ of $C$,
$\Phi(F)$ is a semi-distinguished $A$-invariant face of $K$.}
\par\parindent=16 pt
(b)
\par\parindent=24 pt
(i)
{\it For each nonzero $A$-invariant face $F$ of $K$ associated with $\lambda$,
there \par\parindent=0 pt exists a semi-distinguished $A$-invariant face
$G\subseteq F$ such that $\sp_A(G)=\sp_A(F)$.}
\par\parindent=24 pt
(ii)
{\it For each nonzero $A$-invariant face $F$ of $K$ associated with $\lambda$,
there \par\parindent=0 pt exists in $F$ a generalized eigenvector of $A$
corresponding to $\lambda$
of $\order \nu_\lambda(\left. A\right|_{\span F})$.}
\par\parindent=24 pt
(iii)
{\it For each nonzero $A$-invariant face $F$ of $K$ associated with $\lambda$,
$F$ is a semi-distinguished $A$-invariant face,
if $\sp_A(G)\prec\sp_A(F)$
for all $A$-invariant faces $G$ properly included in $F$.}
\par\parindent=16 pt
(c)
{\it For any $A$-invariant face $F$ of $K$ which contains in its
relative\par\parindent=0 pt
interior a generalized eigenvector of $A$ corresponding
to $\lambda$,
there exists a semi-distinguished $A$-invariant face $G$ included in $F$ such
that $\sp_A(G)=\sp_A(F)$.}
\par\parindent=16 pt
{\it Then conditions} (i),
(ii) {\it given in} (a) ({\it and also conditions} (i)--(iii) {\it given in}
(b))
{\it are equivalent.
Furthermore,
we have} (a)$\Longrightarrow$(c) {\it and} (b)$\Longrightarrow$(c).
\par\parindent=0 pt
\vskip 12 pt

{\it Proof.}
First,
we show that conditions (i),
(ii) of (a) are equivalent.
\par\parindent=16 pt
(i)$\Longrightarrow$(ii):
Let $F$ be a semi-distinguished $B$-invariant face of $C$.
Since $F$ contains in its relative interior a generalized eigen-vector of $A$
corresponding \par\parindent=0 pt to $\lambda$,
so does $\Phi(F)$;
hence,
by condition (i),
we can write $\Phi(F)$ as $G_1\bigvee\cdots\bigvee G_k$,
where each $G_i$ is a semi-distinguished $A$-invariant face of $K$ associated
with $\lambda$.
If $k=1$,
we are done.
So suppose that
$k\ge 2$ and also that each $G_i$ is properly included in $\Phi(F)$.
It is straightforward to show that each
$F_i\bigcap{\cal N}({(\lambda I-A)}^n)$ is a $B$-invariant face of $C$ properly
included in $F$,
and also that the join (in $C$) of these faces is $F$.
This contradicts our assumption on $F$.
\par\parindent=16 pt

(ii)$\Longrightarrow$(i):
Consider any $A$-invariant face $F$ of $K$ which contains in its relative
interior a generalized eigenvector of
$A$ corresponding to $\lambda$.
Let $G$ \par\parindent=0 pt denote $F\bigcap{\cal N}({(\lambda I-A)}^n)$.
Then $G$ is a $B$-invariant face of $C$,
and $\Phi(G)=F$.
Now write $G$ as the join of some $B$-invariant join-irreducible
faces of $C$,
say,
$G_1,\cdots,G_k$.
Since $\lambda$ is the only eigenvalue of $B$,
clearly the faces $G_1,\cdots,G_k$ are all semi-distinguished $B$-invariant.
Then $F=\Phi(G)=\Phi(G_1)\bigvee\cdots\bigvee\Phi(G_k)$,
where each $\Phi(G_i)$ is clearly associated with $\lambda$ and,
by condition (ii),
is a semi-distinguished $A$-invariant face of $K$.
\par\parindent=16 pt

Next,
we show that conditions (i),
(ii) and (iii) of (b) are equivalent.
\par\parindent=16 pt
(i)$\Longrightarrow$(ii):
Let $F$ be an $A$-invariant face associated with $\lambda$,
for which
there exists,
by condition (i),
a semi-distinguished $A$-invariant face $G\subseteq F$ such that
$\sp_A(G)=\sp_A(F)$.
By definition of a semi-distinguished $A$-invariant face,
there exists a generalized eigenvector $x$ of $A$ corresponding to
$\rho_{{}_{\scriptstyle G}}=\rho_{{}_{\scriptstyle F}}=\lambda$ that lies in
$\relint G$.
The order of $x$ as a generalized eigenvector of $A$ is,
by definition of the spectral pair of a face and by Theorem 4.4,
equal to $\nu_\lambda(\left. A\right|_{\span G})$,
and hence to $\nu_\lambda(\left. A\right|_{\span F})$,
as $\sp_A(G)=\sp_A(F)$.
Thus,
$x$ is the desired vector.
\par\parindent=16 pt

(ii)$\Longrightarrow$(iii):
The given spectral pair condition on $F$ clearly implies its $A$-invariant
join-irreducibility.
By condition (ii) $F$ contains at least one \par\parindent=0 pt generalized
eigenvector of $A$,
say $x$,
corresponding to $\rho_{{}_{\scriptstyle F}}$
($=\lambda$)
of $\order\nu_{\rho_{{}_{\scriptscriptstyle F}}}(\left. A\right|_{\span F})$.
Replacing $x$ by ${(I+A)}^{n-1}x$,
if necessary,
we may assume that $\Phi(x)$ is $A$-invariant.
Since $\Phi(x)\subseteq F$ and
$\sp_A(\Phi(x))=(\rho_{{}_{\scriptstyle F}}, \nu_{\rho_{{}_{\scriptscriptstyle
F}}}(\left. A\right|_{\span F}))=\sp_A(F)$,
by the given condition on $F$,
we have $\Phi(x)=F$,
and hence $x\in\relint F$.
This proves that $F$ is semi-distinguished $A$-invariant.
\par\parindent=16 pt

(iii)$\Longrightarrow$(i):
Let $F$ be a nonzero $A$-invariant face associated with $\lambda$.
Choose an $A$-invariant face $G$ included in $F$ minimal with respect to the
property that $\sp_A(G)=\sp_A(F)$.
In view of condition (iii),
$G$ must be semi-distinguished $A$-invariant.
Hence,
condition (i) follows.
\par\parindent=16 pt
(a)(i)$\Longrightarrow$(c):
Follows from Theorem 4.9(i)(c).
\par\parindent=16 pt
(b)(i)$\Longrightarrow$(c):
Obvious.
\par
\hfill{$\Box$}
\par\parindent=16 pt
\vskip 12 pt

The proof of the implication (i)$\Longrightarrow$(ii) in Theorem 6.4 (b)
actually \par\parindent=0 pt establishes the following result:
\par\parindent=0 pt
\vskip 12 pt

{\it Remark} 6.5.$\,\,$
Given a nonzero $A$-invariant face $F$ of $K$,
if there exists a semi-distinguished $A$-invariant face $G$ included in $F$
such that $\sp_A(G)=\sp_A(F)$,
then there exists in $F$ a generalized eigenvector of $A$ corresponding to
$\rho_{{}_{\scriptstyle F}}$ of
$\order\,\nu_{\rho_{{}_{\scriptscriptstyle F}}}(\left. A\right|_{\span F})$.
\par\parindent=16 pt
\vskip 12 pt

As we shall show in Example 8.1 of Section 8,
the converse of Remark 6.5 does not hold.
\par\parindent=16 pt

Some further remarks about the conditions (a),
(b) and (c) of Theorem 6.4 are in order.
In general,
we do not have the implication (a)$\Longrightarrow$(b).
As one can readily check,
for the $A$ and $K$ of Example 5.3,
with $\lambda=\rho(A)$,
condition (a)(i) is satisfied,
but not condition (b)(i).
Whether the implication (b)$\Longrightarrow$(a) holds or not is {\it unknown}
to us.
In Example 8.2 we shall show that the implications (c)$\Longrightarrow$(a)
and (c)$\Longrightarrow$(b) both do not hold.
\par\parindent=16 pt

If $K$ is a proper cone with the property that the dual cone of each of its
faces
is a facially exposed cone and if $A\in\pi(K)$,
then according to Theorem 6.1,
for any nonzero $A$-invariant face $F$ of $K$,
$F$ is semi-distinguished $A$-invariant
if and only if $\sp_A(F)\succ\sp_A(G)$ for all $A$-invariant faces $G$
properly included in $F$;
thus,
for all distinguished eigenvalues $\lambda$ of $A$,
condition (b)
(and hence also condition (c))
of Theorem 6.4 is satisfied.
Clearly proper cones $K$ all of whose nontrivial faces are polyhedral
and for which $K^*$ are facially exposed cones have the property.
Such cones include all polyhedral cones,
proper cones of $\IR^{\,\,3}$ whose dual cones are facially exposed,
and also strictly convex smooth cones.
[Recall that a proper cone $K$ is said to be {\it strictly convex} if each of
its
nontrivial faces is one-dimensional,
and is {\it smooth} if for each nonzero boundary vector $x$ of $K$,
$\dim d_K(\Phi(x))=1$.
It can be shown that the class of strictly convex cones whose dual cones
are facially exposed is the same as the class of strictly convex smooth
cones.].
\par\parindent=16 pt

Perfect cones form another interesting class with the above-mentioned property.
Following [Bar 1],
we call a proper cone a {\it perfect cone} if each of its faces is self-dual in
its own linear span.
Examples of perfect cones include,
the nonnegative orthant $\IR^{\,\,n}_{\,\,+}$,
the $n$-dimensional ice-cream cone $K_n:=\{ {(\xi_1,\cdots,\xi_n)}^T: \xi_1\ge
{({\xi_2}^2+\cdots+{\xi_n}^2)}^{1/2}\}$,
and the cone of $n\!\times\! n$ positive semi-definite hermitian (or real
symmetric) matrices.
An equivalent definition for $K$ to be a perfect cone is that $K$ is self-dual,
and for each face $F$ of $K$,
the orthogonal projection $p_{{}_{\scriptstyle F}}$ of $\IR^{\,\,n}$ onto
$\span F$ belongs to $\pi(K)$
(see [Bar 1, Proposition 1]).
So a perfect cone is an orthogonal projectionally exposed,
and hence a projectionally exposed
cone (see [Su--T] for definitions).
But by [Su--T, Corollary 4.4 and Lemma 2.2] each face of a projectionally
exposed cone is a facially
exposed cone in its own right.
Hence,
each face of a perfect cone is itself a facially exposed self-dual cone.
\par\parindent=16 pt

For each nonnegative integer $n$,
we use $P(n)$ to denote the cone of all real polynomials of degree not
exceeding
$n$ that are nonnegative on the closed interval $[0,1]$.
In [B--T] it is proved that ${P(n)}^*$ is always a facially exposed cone,
and also that each nonzero face of $P(n)$ is linearly isomorphic with $P(m)$
for some nonnegative integer $m\le n$.
Thus,
the cone $P(n)$ also serves as another interesting example of a proper cone
with the property that the dual cone of each of its faces is a facially exposed
cone.
\par\parindent=16 pt

To sum up our preceding discussions,
we have the following result.
\par\parindent=0 pt
\vskip 12 pt

{\bf Theorem 6.6.}
\par\parindent=16 pt
(i)
{\it Let $K$ be a proper cone with the property that the dual cone of
each of its faces is a facially exposed cone,
and let $A\in\pi(K)$.
Then for any nonzero $A$-invariant face $F$ of $K$,
$F$ is semi-distinguished $A$-invariant if and only if
$\sp_A(G)\prec\sp_A(F)$ for all $A$-invariant faces $G$ properly included
in $F$.
Consequently,
for all distinguished eigenvalues $\lambda$ of $A$,
condition} (b) ({\it and hence} (c))
{\it of Theorem 6.4 is satisfied.}

\par\parindent=16 pt
(ii)
{\it A proper cone $K$ has the property given in the hypothesis of part} (i),
{\it if it fulfills one of the following conditions:}
\par\parindent=24 pt
(a) {\it $K^*$ is a facially exposed cone,
and all nontrivial faces of $K$ are polyhedral
(which is the case,
if $K$ is polyhedral,
or is a strictly convex smooth cone).}
\par\parindent=24 pt
(b) {\it $K$ is a perfect cone.}
\par\parindent=24 pt
(c) {\it $K$ equals $P(n)$ for some nonnegative integer $n$.}
\par\parindent=16 pt
\vskip 12 pt

Note that we could have formulated the last part of Theorem 6.6
(i)\par\parindent=0 pt
without mentioning the distinguished eigenvalues of $A$.
Indeed,
the statement that,
for all distinguished eigenvalues $\lambda$ of $A$,
condition (b) (ii) of Theorem 6.4 is satisfied amounts to saying that for each
nonzero $A$-invariant face $F$ of $K$,
there exists in $F$ a distinguished generalized eigenvector of $A$
corresponding to $\rho_{{}_{\scriptstyle F}}$ of
$\order\nu_{\rho_{{}_{\scriptscriptstyle F}}}(\left. A\right|_{\span F})$.
We leave it to the reader to give similar \par reformulations for the statement
that conditions (b) (ii),
(iii) or (c) of \par Theorem 6.4 are satisfied for all distinguished
eigenvalues
$\lambda$ of $A$.
\par\parindent=0 pt
\vskip 12 pt

{\bf Theorem 6.7.}
{\it Let $A\in\pi(K)$,
where $K$ is a polyhedral cone or a strictly convex cone.
Then for any distinguished eigenvalue $\lambda$ of $A$,
condition (a) of Theorem 6.4 is satisfied.}
\par\parindent=16 pt
\vskip 12 pt

The polyhedral case of Theorem 6.7 is covered by Lemma 6.8.
The proof for a strictly convex cone will follow from Lemma 6.10.
\par\parindent=0 pt
\vskip 12 pt

{\bf Lemma 6.8.}
{\it Let $K$ be a polyhedral cone.
Let $A\in\pi(K)$,
and let $\lambda$ be a \par\parindent=0 pt distinguished eigenvalue of $A$ for
$K$.
Let $C$ denote the cone $K\bigcap{\cal N}({(\lambda I-A)}^n)$,
and let $B$ denote the restriction map $\left. A\right|_{\span C}$.
Then for any semi-distinguished $B$-invariant face $F$ of $C$,
$\Phi(F)$ is a semi-distinguished $A$-invariant face of $K$.}
\par\parindent=0 pt
\vskip 12 pt

{\it Proof.}
Since $K$ is polyhedral,
by Theorem 6.6 it suffices to show that $\sp_A(G)\prec\sp_A(\Phi(F))$
for any $A$-invariant face $G$ properly included in $\Phi(F)$.
If $\rho_{{}_{\scriptstyle G}}<\lambda$,
clearly there is no problem.
So,
suppose that $\rho_{{}_{\scriptstyle G}}=\lambda$.
By Theorem 6.6 again,
condition (b) (i) of Theorem 6.4 is satisfied.
Hence,
there exists a semi-distinguished $A$-invariant face $H$
included in $G$ such that $\sp_A(H)=\sp_A(G)$.
Then $H\bigcap{\cal N}({(\lambda I-A)}^n)$ is a $B$-invariant face
of $C$ \par\parindent=0 pt properly included in $F$.
Since $F$ is a semi-distinguished $B$-invariant face of the polyhedral cone
$C$,
by Theorem 6.6 again we have $\sp_B(H\bigcap{\cal N}({(\lambda I-A)}^n))$
$\prec\sp_B(F)$.
Now $\relint F\subseteq\relint\Phi(F)$ and $B$ is a restriction of $A$;
so by the definition of the spectral pair of a face,
we have $\sp_B(F)=\sp_A(\Phi(F))$.
But $H$ contains in its relative interior a generalized eigenvector of $A$
corresponding to $\lambda$,
so for a similar reason,
we also have $\sp_B(H\bigcap{\cal N}({(\lambda I-A)}^n))=\sp_A(H)$.
Hence,
the desired inequality $\sp_A(G)\prec\sp_A(\Phi(F))$ follows.
\par
\hfill{$\Box$}
\par\parindent=16 pt
\vskip 12 pt

For completeness,
we also note that the converse of Lemma 6.8 always holds,
even without the polyhedrality assumption on $K$.
The proof is straightforward.
\par\parindent=0 pt
\vskip 12 pt

{\it Remark} 6.9.$\,\,$
Let $A\in\pi(K)$ and let $\lambda$ be a distinguished eigenvalue of $A$ for
$K$.
Let $C$ denote the cone $K\bigcap{\cal N}({(\lambda I-A)}^n)$ and let $B$
denote the restriction map $\left. A\right|_{\span C}$.
For any face $F$ of $C$,
if $\Phi(F)$ is a semi-distinguished $A$-invariant face of $K$,
then $F$ is a semi-distinguished $B$-invariant face of $C$.
\par\parindent=0 pt
\vskip 12 pt

{\bf Lemma 6.10.}
{\it Let $K$ be a strictly convex cone,
and let $A\in\pi(K)$.
For any $A$-invariant face $F$ of $K$,
we have $\sp_A(F)\succ\sp_A(G)$ for all $A$-invariant faces $G$ properly
included in $F$,
whenever $F$ is a semi-distinguished $A$-invariant face of $K$.}
\par\parindent=0 pt
\vskip 12 pt

{\it Proof.}
It is clear that each nontrivial $A$-invariant face of $K$,
being an extreme ray of $K$,
is a distinguished $A$-invariant face and satisfies trivially the desired
property involving spectral pairs.
So it suffices to consider the case $F=K$.
If $\nu_{\rho(A)}(A)>1$,
then for any nonzero $A$-invariant face $G$ properly included in $K$,
we have,
\begin{center}
$\sp_A(K)=(\rho(A),\nu_{\rho(A)}(A))\succ (\rho(A),1)\succeq
(\rho_{{}_{\scriptstyle G}},1)=\sp_A(G)$,
\end{center}
\par\parindent=0 pt

where the last equality holds as $G$ is an extreme ray.
So suppose that $\nu_{\rho(A)}(A)=1$.
Then the generalized eigenvector of $A$
(corresponding to $\rho(A)$)
that lies in $\int K$ is in fact an eigenvector.
Hence,
condition (d) of Theorem 4.11 is satisfied;
that is,
$K$ is a distinguished $A$-invariant face of itself,
and so we also have $\sp_A(K)\succ\sp_A(G)$ for any $A$-invariant
face $G$ properly included in $K$.
\par
\hfill{$\Box$}
\par\parindent=16 pt
\vskip 12 pt

We now give the proof of Theorem 6.7 for the strictly convex case:
\par\parindent=16 pt
\vskip 12 pt

Clearly it suffices to show that if $K$ contains in its interior
a generalized eigenvector of $A$ then $K$ can be expressed as a join of
semi-distinguished $A$-invariant faces associated with $\rho(A)$.
If $\nu_{\rho(A)}(A)>1$ then,
as shown in the proof of Lemma 6.10,
we have $\sp_A(K)\succ\sp_A(G)$ for all $A$-invariant faces $G$ properly
included in $K$;
hence,
$K$ is $A$-invariant join-irreducible,
and since $K$ contains in its interior a generalized eigenvector of $A$,
$K$ is itself a semi-distinguished $A$-invariant face.
When $\nu_{\rho(A)}(A)=1$,
we have either $K$ is a distinguished $A$-invariant face of itself,
or $\dim (K\bigcap{\cal N}(\rho(A)I-A))\ge 2$.
In the former case,
we are done.
In the latter case,
since $K$ is strictly convex,
clearly $K$ can be expressed as the join of two $A$-invariant extreme rays,
and hence of two $A$-invariant distinguished faces,
associated with $\rho(A)$.
\par\parindent=16 pt
\vskip 12 pt

For completeness,
we also include the following interesting result on a smooth cone.
\par\parindent=0 pt
\vskip 12 pt

{\bf Lemma 6.11.}
{\it Let $K$ be a smooth cone,
and let $A\in\pi(K)$.}
\par\parindent=16 pt
(i)
{\it If $\nu_{\rho(A)}(A)>1$,
then $K$ contains in its interior a generalized \par\parindent=0 pt eigenvector
of $A$.}
\par\parindent=16 pt
(ii)
{\it There always exists in $K$ a generalized eigenvector of $A$ of
$\order\nu_{\rho(A)}(A)$ corresponding to $\rho(A)$.}
\par\parindent=16 pt
(iii)
{\it If $\sp_A(K)\succ\sp_A(G)$ for all $A$-invariant faces $G$ properly
included in $K$,
then $K$ is a semi-distinguished $A$-invariant face of itself.}
\par\parindent=0 pt
\vskip 12 pt

{\it Proof.}
\par\parindent=16 pt
(i):
Assume to the contrary that
$K$ does not contain in
its interior a \par\parindent=0 pt generalized eigenvector of $A$.
By [Tam 3,
Theorem 5.1]
$A^T$ has a distinguished eigenvalue $\lambda$ other than $\rho(A)$.
Let $w\in K^*$ be a corresponding eigenvector.
Choose from $K$ an eigenvector $x$ of $A$ corresponding to $\rho(A)$,
which is \par associated with a maximal Jordan block
(see Remark 4.2).
Also choose from $K^*$ an eigenvector $z$ of $A^T$ corresponding to $\rho(A)$,
which is again associated with a maximal Jordan block.
Since $x$ and $w$ are respectively eigenvectors of $A$ and $A^T$ corresponding
to distinct eigenvalues,
$x$ and $w$ are mutually orthogonal.
Because the size of a maximal Jordan block of $A$ corresponding to $\rho(A)$,
i.e.
$\nu_{\rho(A)}(A)$,
is greater than $1$,
the vectors $x$ and $z$ are also mutually orthogonal.
Thus $x$ is orthogonal to two linearly independent vectors of $K^*$.
This contradicts the smoothness of $K$ at $x$.
\par\parindent=16 pt
(ii):
If $\nu_{\rho(A)}(A)=1$,
clearly there is no problem.
If $\nu_{\rho(A)}(A)>1$,
our assertion follows from part (i).
\par\parindent= 16 pt
(iii):
If $\nu_{\rho(A)}(A)=1$,
then by the spectral pair condition on $K$,
clearly $K$ is a distinguished $A$-invariant face of itself.
Otherwise,
by part (i) $K$ contains in its interior a generalized
eigenvector of $A$.
And since the spectral pair condition on $K$ implies the $A$-invariant
join-irreducibility of $K$,
it follows that $K$ is a semi-distinguished $A$-invariant face of itself.
\par
\hfill{$\Box$}
\par\parindent=16 pt
\vskip 12 pt

The following example shows that parts (i)--(iii) of Lemma 6.11 all do not
hold,
if in the hypothesis on $K$,
$\lq\lq$smooth$"$ is replaced by $\lq\lq$strictly convex$"$.
Consequently,
the converse of Lemma 6.10 also does not hold.
\par\parindent=0 pt
\vskip 12 pt

{\bf Example 6.12.}$\,\,$
Let $\alpha$ be a positive real number less than $1$.
Let $f$ be a real-valued function defined on the real line that possesses
the following \par\parindent=0 pt properties:
\par\parindent=16 pt
(i)
$\,\,\,f$ is a differentiable convex function;
\par\parindent=16 pt
(ii)
${{}\atop{{\textstyle\lim}\atop{\scriptstyle\xi\to -\infty}}}f(\xi)=\infty$;
\par\parindent=16 pt
(iii)
${{}\atop{{\textstyle\lim}\atop{\scriptstyle\xi\to \infty}}}f(\xi)=-\infty$;
\par\parindent=16 pt
(iv)
${{}\atop{{\textstyle\lim}\atop{\scriptstyle\xi\to \infty}}}f^\prime(\xi)=0$;
$\,\,\,$and
\par\parindent=16 pt
(v)$\,\,\,$
For each $\xi\in\IR\,\,$,
we have $f(\xi+1)\le\alpha f(\xi)$.
\par\parindent=16 pt

For instance,
if $e^{-1}<\alpha<1$ and $f$ is given by:
\begin{center}
$f(\xi)=
\left\{
\right.
\begin{array}{rcl}
-\log (\xi+1) & \mbox{for} & \xi\ge 0 \cr
\exp (-\xi)-1 & \mbox{for} & \xi<0
\end{array}
$
,
\end{center}
\par\parindent=0 pt
then one can readily verify that $f$ possesses properties (i)--(v).
\par\parindent=16 pt

Denote the epigraph of $f$,
i.e.
the set $\{ {(\xi_1,\xi_2)}^T\in\IR^{\,\,2}: \xi_2\ge f(\xi_1)\}$,
by $\epi f$.
Since $f$ is a continuous convex function,
$\epi f$ is a closed convex set.
By the convexity of $f$ and property (v),
for any $\xi\in\IR\,\,$,
we have,
$f^\prime(\xi)\le f(\xi+1)-f(\xi)\le -(1-\alpha)f(\xi)$.
In view of property (ii) it follows that
$\lim_{\,\xi\to -\infty}f^\prime(\xi)=-\infty$.
The latter condition implies that the recession cone ${\bf 0}^+(\epi f)$
cannot contain a vector of the form ${(\nu_1,\nu_2)}^T$ with $\nu_1<0$.
Similarly,
property (iv) implies that ${\bf 0}^+(\epi f)$
cannot contain a vector of the form ${(\nu_1,\nu_2)}^T$ with
$\nu_2<0$.
Clearly,
${(0,1)}^T\in{\bf 0}^+(\epi f)$.
Also,
we have,
${(1,0)}^T\in{\bf 0}^+(\epi f)$;
otherwise,
$f$ is not one-to-one and property (iii) cannot be met.
Thus,
${\bf 0}^+(\epi f)$ equals $\,\IR^{\,\,2}_{\,\,+}$.
\par\parindent=16 pt

Let $K$ be the proper cone of $\IR^{\,\,3}$ whose dual cone $K^*$ is given by:
\begin{center}
$K^*=\{\lambda{(\xi_1,\xi_2,1)}^T\in\IR^{\,\,3}: \lambda\ge 0,
\,{(\xi_1,\xi_2)}^T\in\epi(f)\}\,\bigcup\, ({\bf 0}^+(\epi f)\times\{ 0\})$.
\end{center}
\par\parindent=0 pt

Since $f$ is a differentiable function,
$K^*$ is smooth at each nonzero boundary vector of the form
$\lambda{(\xi_1,\xi_2,1)}^T$.
To show that $K^*$ is smooth at ${(1,0,0)}^T$,
let ${(\eta_1,\eta_2,\eta_3)}^T$ be a nonzero vector of $K$ which is orthogonal
to
${(1,0,0)}^T$.
Clearly,
$\eta_1=0$.
Since ${(0,1,0)}^T\in K^*$,
$\eta_2\ge 0$.
If $\eta_2>0$,
then in view of property (iii),
we can find some $\xi\in\IR\,$ such that the inner product between
${(\eta_1,\eta_2,\eta_3)}^T$ and ${(\xi,f(\xi),1)}^T$ is negative,
which is a contradiction.
It follows that
(up to multiples)
${(0,0,1)}^T$ is the only vector of $K$ orthogonal to ${(1,0,0)}^T$.
Similarly,
we can show that ${(0,0,1)}^T$ is also the only vector of $K$ orthogonal to
${(0,1,0)}^T$.
This shows that $K^*$ is a smooth cone.
Hence,
$K$ is a strictly convex cone.
\par\parindent=16 pt

Now let
\begin{center}
$A=
\left[
\begin{array}{ccc}
1 & 0	   & 0 \\
0 & \alpha & 0 \\
1 & 0	   & 1
\end{array}
\right]
$.
\end{center}
\par\parindent=0 pt
By property (v) one readily checks that $A^T\in\pi(K^*)$;
hence,
$A\in\pi(K)$.
Clearly,
$\nu_{\rho(A)}(A)=2$.
Note that $K$ does not contain a generalized eigenvector of $A$ corresponding
to $\rho(A)$ of order two;
because,
any such vector is necessarily orthogonal to the vector ${(0,1,0)}^T$ of $K^*$,
contradicting what we have found above.
It is now clear that $K$ is not a semi-distinguished $A$-invariant face
of itself,
and also that we have $\sp_A(K)\succ\sp_A(G)$ for all $A$-invariant faces $G$
properly included in $K$.
\par\parindent=0 pt
\vskip 25 pt

\begin{center}
7. E{\scriptsize XTENSIONS OF

TWO NONNEGATIVE MATRIX RESULTS}
\end{center}
\par\parindent=16 pt

Besides Theorem B (mentioned in Section 1),
the following is another early result on the combinatorial spectral theory
of a nonnegative matrix:
\par\parindent=0 pt
\vskip 12 pt

{\bf Theorem C.} ([Schn 1, Theorem 3]).
{\it If $P$ is a nonnegative matrix,
then $\nu_{\rho(P)}(P)=1$ if and only if there are no comparable basic classes
of $P$.}
\par\parindent=16 pt
\vskip 12 pt

It is clear that Theorem C follows as a special case of Theorem A,
the Rothblum Index theorem.
\par\parindent=16 pt

For a general cone-preserving map,
we have the following result related to Theorem C.
\par\parindent=0 pt
\vskip 12 pt

{\bf Theorem 7.1.}
{\it Let $A\in\pi(K)$,
and let $\lambda$ be a distinguished eigenvalue of $A$ for $K$.
Denote by $m_\lambda$ the maximal order of distinguished generalized
eigenvectors of $A$ corresponding to $\lambda$.
Consider the following
conditions:} \par\parindent=16 pt
(a) {\it $m_\lambda=1$.}
\par\parindent=16 pt
(b) {\it Any two distinct semi-distinguished $A$-invariant faces of $K$
associated
with $\lambda$ are noncomparable.}
\par\parindent=16 pt
(i) {\it $\,$We always have} (a)$\Longrightarrow$(b).
\par\parindent=16 pt
(ii) {\it When condition} (c)
({\it and hence also when conditions} (a) {\it or} (b))
{\it of Theorem 6.4 is satisfied,
conditions} (a) {\it and} (b) {\it are equivalent.}
\par\parindent=0 pt
\vskip 12 pt

{\it Proof.}
\par\parindent=16 pt
(i):
It is clear that any two distinct distinguished $A$-invariant faces of $K$
associated with $\lambda$ are noncomparable.
So it suffices to show that,
when $m_\lambda=1$,
any semi-distinguished $A$-invariant face associated with $\lambda$ is
distinguished $A$-invariant.
Indeed,
if $F$ is a semi-distinguished $A$-invariant face
associated with $\lambda$,
then by definition $F$ is $A$-invariant join-irreducible and contains in its
relative interior a generalized eigenvector
of $A$ corresponding to $\lambda$.
Since $m_\lambda=1$,
this generalized eigenvector
is necessarily an eigenvector.
Hence,
by Theorem 4.11 $F$ must be a distinguished $A$-invariant face of $K$
associated with $\lambda$.
This establishes (b).
\par\parindent=16 pt

(ii):
It suffices to establish the implication (b)$\Longrightarrow$(a).
Assume to the \par\parindent=0 pt contrary that condition (b) is satisfied,
but not condition (a).
Then $m_\lambda>1$.
Let $x$ be a distinguished generalized eigenvector of $A$ corresponding to
$\lambda$ of $\order m_\lambda$.
Then $F=\Phi({(I+A)}^{n-1}x)$ is an $A$-invariant face that contains in \par
its relative
interior a generalized eigenvector of $A$ corresponding to $\lambda$.\par
Furthermore,
$\sp_A(F)=(\lambda, m_\lambda)$.
By condition (c) of Theorem 6.4 there exists a semi-distinguished $A$-invariant
face $G$ included in $F$ such that $\sp_A(G)=\sp_A(F)$.
Now choose any extremal distinguished eigenvector $u$ of $A$ \par corresponding
to $\lambda$ that lies in $G$.
Then by Theorem 4.11 $\Phi(u)$ is a \par distinguished,
and hence semi-distinguished,
$A$-invariant face associated with $\lambda$ which is properly included in $G$.
This contradicts condition (b).
\par
\hfill{$\Box$}
\par\parindent=16 pt
\vskip 12 pt

Later,
by Example 8.1
we shall show that for the conditions (a),
(b) of Theorem 7.1,
the implication (b)$\Longrightarrow$(a) is not valid
in general.

\par\parindent=16 pt
In related to Theorem B we have the following result.
\par\parindent=0 pt
\vskip 12 pt

{\bf Theorem 7.2.}
{\it Let $A\in\pi(K)$,
and let $\lambda$ be a distinguished eigenvalue of $A$ for $K$.
Consider the following conditions:}
\par\parindent=16 pt
(a) {\it $\dim{\cal N}(\lambda I-A)=1$.}
\par\parindent=16 pt
(b) {\it Any $A$-invariant face which contains in its relative interior a
\par\parindent=0 pt generalized eigenvector of $A$ corresponding to $\lambda$
is semi-distinguished
\par
$A$-invariant.}
\par\parindent=16 pt
(c) {\it Any two $A$-invariant faces of $K$ each containing in its relative
interior
a generalized eigenvector of $A$ corresponding to $\lambda$ are comparable.}
\par\parindent=16 pt
(d) {\it Any two semi-distinguished $A$-invariant faces of $K$ associated with
$\lambda$
are comparable.}
\par\parindent=16 pt
(e) {\it $\dim\span [{\cal N}(\lambda I-A)\bigcap K]=1$.}
\par\parindent=16 pt
(i) {\it We always have the implications}
\par\parindent=35 pt
(a)$\Longrightarrow$(b)$\Longrightarrow$(c)$\Longrightarrow$(d)$\Longrightarrow
$(e).
\par\parindent=16 pt
(ii) {\it When condition} (a) {\it of Theorem 6.4 is satisfied,
conditions} (b),
(c) {\it and} (d) {\it are equivalent.}
\par\parindent=16 pt
(iii) {\it When $K$ is polyhedral,
conditions} (b),
(c) {\it and} (d) {\it are equivalent.
If,
in addition,
$\lambda=\rho(A)$,
then condition} (a) {\it is also another equivalent condition.}
\par\parindent=0 pt
\vskip 12 pt

{\it Proof.}
\par\parindent=16 pt
(i):
\par\parindent=24 pt
(a)$\Longrightarrow$(b):
We are required to show that if $F$ is any $A$-invariant face of $K$ that
contains in its relative interior a generalized eigenvector of $A$
corresponding to $\lambda$,
then $F$ is $A$-invariant join-irreducible.
By considering the restriction map $\left. A\right|_{\span F}$ instead of $A$,
we may henceforth take $F$ to be $K$ and $\lambda$ to be $\rho(A)$.
Let $G$ be any $A$-invariant face of $K$ properly included in $K$.
We are going to show that $\sp_A(G)\prec\sp_A(K)$.
Once this is done,
the $A$-invariant join-irreducibility of $K$ will follow.
\par\parindent=16 pt

Now $\span G$ is an $A$-invariant subspace of $\IR^{\,\,n}$.
If $\span G$
contains a \par\parindent=0 pt generalized eigenvector of $A$ corresponding to
$\rho(A)$ of $\order\nu_{\rho(A)}(A)$,
then $\span G$ will include ${\cal N}({(\rho(A)I-A)}^n)$,
as the Jordan form of $A$ has only one block corresponding to $\rho(A)$.
Since $G=K\bigcap\span G$ and ${\cal N}({(\rho(A)I-A)}^n)$ meets the interior
of $K$,
it will follow that $G$ contains an interior vector of $K$ and hence is equal
to $K$,
which is a contradiction.
We just show that if $\rho(A)$ is an eigenvalue of $\left. A\right|_{\span G}$
then its index
is less than $\nu_{\rho(A)}(A)$.
Now it is clear that we have $\sp_A(G)\prec\sp_A(K)$.
\par\parindent=24 pt

(b)$\Longrightarrow$(c):
Assume to the contrary that there exist two noncomparable \par\parindent=0 pt
$A$-invariant faces of $K$,
say $F_1, F_2$,
each containing in its relative interior a generalized eigenvector of $A$
corresponding to $\lambda$.
Then $F_1\bigvee F_2$ is an $A$-invariant face with the same property,
but is not semi-distinguished $A$-invariant.
This contradicts condition (b).
\par\parindent=24 pt

(c)$\Longrightarrow$(d):
Obvious.
\par\parindent=24 pt

(d)$\Longrightarrow$(e):
If $\dim [{\cal N}(\lambda I-A)\bigcap K]\ge 2$,
choose any two distinct \par\parindent=0 pt extreme vectors $x_1, x_2$ of the
cone
${\cal N}(\lambda I-A)\bigcap K$.
Then $\Phi(x_1)$ and $\Phi(x_2)$ are distinguished,
and hence semi-distinguished,
noncomparable $A$-invariant faces of $K$ associated with $\lambda$.
This violates condition (d).
\par\parindent=16 pt

(ii):
To prove the equivalence of conditions (b),
(c),
(d),
it suffices to establish the implication (d)$\Longrightarrow$(b).
Let $F$ be an $A$-invariant face which contains in its relative interior a
generalized eigenvector of $A$ corresponding to $\lambda$.
By condition (a)(i) of Theorem 6.4 we can write $F$ as $F_1\bigvee\cdots\bigvee
F_k$,
where each $F_i$ is a semi-distinguished $A$-invariant face
associated with $\lambda$.
We may assume that
there is no redundant terms in this representation of $F$.
If $k\ge 2$,
then $F_1$ and $F_2$ will be two noncomparable semi-distinguished $A$-invariant
faces
associated with $\lambda$,
contradicting condition (d).
Hence,
we have $k=1$ and $F$ is semi-distinguished $A$-invariant.
\par\parindent=16 pt

(iii):
When $K$ is polyhedral,
by Theorem 6.7 condition (a) of Theorem 6.4 is satisfied.
Hence,
by part (ii),
conditions (b),
(c) and (d) are equivalent.
\par\parindent=16 pt

Now suppose,
in addition,
that $\lambda=\rho(A)$.
To prove the equivalence of conditions (a)--(d),
it remains to establish the implication
(d)$\Longrightarrow$(a). \par\parindent=16 pt

Suppose that condition (d) is satisfied.
Clearly an equivalent condition \par\parindent=0 pt for (a) is that the
geometric multiplicity of
$\rho(A)$ as an eigenvalue of $\left. A\right|_{{\cal N}({(\rho(A)I-A)}^n)}$
$\,\,$ is $\,$ one.
Since $K$ is polyhedral,
${\cal N}({(\rho(A)I-A)}^n)$ has a basis that consists \par of vectors of $K$
(see [Tam 3,
Theorem 7.5]);
so
${\cal N}({(\rho(A)I-A)}^n)\bigcap K$ is a proper cone of
${\cal N}({(\rho(A)I-A)}^n)$.
Denote by $B$ the restriction of $A$ to ${\cal N}({(\rho(A)I-A)}^n)$ and by
$C$ the cone ${\cal N}({(\rho(A)I-A)}^n)\bigcap K$.
If $G_1, G_2$ are semi-distinguished $B$-invariant faces of $C$,
then by Lemma 6.8 $\Phi(G_1)$ and $\Phi(G_2)$ are semi-distinguished
$A$-invariant faces of $K$ and are comparable \par\parindent=0 pt according to
condition (d).
But $G_i={\cal N}({(\rho(A)I-A)}^n)\bigcap\Phi(G_i)$ for $i=1, 2$,
hence $G_1$ and $G_2$ are also comparable.
This shows that $B$
[$\in\pi(C)$] also satisfies the condition
corresponding to condition (d).
Therefore,
by considering $B$ instead of $A$,
we may,
henceforth,
assume that $\rho(A)$ is the only eigenvalue of $A$.
\par\parindent=16 pt

Since $K$ is polyhedral,
by Theorem 6.7,
condition (a)(i) of Theorem 6.4 \par\parindent=0 pt is satisfied;
hence,
any nonzero
$A$-invariant face of $K$ (which necessarily \par contains in its
relative interior a generalized eigenvector of $A$ corresponding to $\rho(A)$)
is the join of semi-distinguished $A$-invariant faces of $K$
(associated with $\rho(A)$).
But according to condition (d),
any two semi-distinguished $A$-invariant faces of $K$ are comparable.
It follows that each nonzero $A$-invariant face of $K$ is semi-distinguished.
In particular,
$K$ is a semi-distinguished $A$-invariant face of itself.
Hence,
by Theorem 6.1 $A^T$ has
(up to multiples)
a unique distinguished eigenvector for $K^*$,
say $w$.
Indeed,
as noted in the proof of Lemma 5.2,
in this case $A^T$ has precisely one maximal Jordan block \par corresponding to
$\rho(A)$ and $w$ is associated with this maximal block.
Furthermore,
since $\rho(A)$ is the only eigenvalue of $A$,
we have ${(\span\{ w\})}^{\bot}=E_{\nu-1}$,
where $E_{\nu-1}$ is the space of generalized eigenvectors of $A$
(corresponding to $\rho(A)$) of order less than or equal
to $\nu_{\rho(A)}(A)-1$,
together with the zero vector.
We contend that
$\Phi(w)$ is an extreme ray of $K^*$.
\par\parindent=16 pt

To show this,
choose a nonzero extreme vector $z$ of $\Phi(w)$.
Then we have
\begin{center}
$(\rho(A),1)=\sp_{A^T}(\Phi(w))\succeq\sp_{A^T}(\Phi(z))\succeq(\rho(A),1)$,
\end{center}
\par\parindent=0 pt
where the last inequality follows from the assumption that $\rho(A)$ is the
only
eigenvalue of $A$.
Hence,
$\sp_{A^T}(\Phi(z))=(\rho(A),1)$ and so $z$ is an eigenvector of $A^T$.
It follows that $z$ is a multiple of $w$.
This proves our contention.
\par\parindent=16 pt

Since $K$ is polyhedral and $\Phi(w)$ is an extreme ray of $K^*$,
$d_{K^*}(\Phi(w))$ must be a cone of dimension $n-1$
(for reference,
see [Tam 1,
Theorem 3]),
and hence is a proper cone of $E_{\nu-1}$;
thus we can write $\left. A\right|_{E_{\nu-1}}\in\pi(d_{K^*}(\Phi(w)))$.
Now let the Jordan form of $A$ be
\begin{center}
$J_{i_1}(\rho(A))\bigoplus J_{i_2}(\rho(A))\bigoplus\cdots\bigoplus
J_{i_k}(\rho(A))$,
\end{center}
\par\parindent=0 pt
where $k\ge 1$,
and $i_1> i_2\ge\cdots\ge i_k$ if $k\ge 2$.
Then it is not difficult to show that the Jordan form of $\left.
A\right|_{E_{\nu-1}}$
is
\begin{center}
$J_{i_1-1}(\rho(A))\bigoplus J_{i_2}(\rho(A))\bigoplus
J_{i_3}(\rho(A))\bigoplus\cdots\bigoplus J_{i_k}(\rho(A))$.
\end{center}
\par\parindent=0 pt
Note that $d_{K^*}(\Phi(w))$ is itself a polyhedral cone,
and also that $\left. A\right|_{E_{\nu-1}}$ inherits from $A$ the condition
corresponding
to condition (d).
So,
applying our above argument with $\left. A\right|_{E_{\nu-1}}$ in place of
$A$,
we see that the Jordan form of $\left. A\right|_{E_{\nu-1}}$ also has precisely
one maximal block.
Suppose $k\ge 2$.
If $i_1-1=i_2$,
we already obtain a contradiction.
Otherwise,
our argument shows that $\left. A\right|_{E_{\nu-2}}$,
where $E_{\nu-2}$ denotes the space of generalized eigenvectors of $A$ of
order less than or equal to $\nu_{\rho(A)}(A)-2$,
together with the zero vector,
belongs to $\pi(G)$ for some $A$-invariant face $G$ of $K$ of dimension $n-2$.
Furthermore,
$\left. A\right|_{E_{\nu-2}}$ also has precisely one maximal block in its
Jordan form,
which is
\begin{center}
$J_{i_1-2}(\rho(A))\bigoplus J_{i_2}(\rho(A))\bigoplus\cdots\bigoplus
J_{i_k}(\rho(A))$.
\end{center}
\par\parindent=0 pt
If $i_1-2=i_2$,
we again arrive at a contradiction.
Otherwise,
we repeat the process and obtain a contradiction after a finite number of
steps.
This shows that $k=1$,
i.e. $A$ has precisely one Jordan block (corresponding to $\rho(A)$),
as desired.
\par
\hfill{$\Box$}
\par\parindent=16 pt
\vskip 12 pt

In related to the conditions of Theorem 7.2 one may also consider the
following condition:
\par\parindent=16 pt

(f) Any two $A$-invariant faces of $K$ associated with $\lambda$  are
comparable.
\par\parindent=16 pt

Clearly condition (f) implies condition (c) of Theorem 7.2.
On the other hand,
condition (f) is logically independent of condition (a) of Theorem 7.2.
To see why condition (a) does not imply condition (f),
consider the $3\!\times\! 3$ nonnegative matrix $A$ given by:
\begin{center}
$A=
\left[
\begin{array}{ccc}
1 & 1 & 1 \\
0 & 0 & 0 \\
0 & 0 & 0
\end{array}
\right]
$
.
\end{center}
\par\parindent=0 pt
Take $\lambda=\rho(A)$.
Clearly $\rho(A)=1$ and the geometric multiplicity of $1$ as an
eigenvalue of $A$ is $1$.
However,
$\Phi(e_1+e_2)$ and $\Phi(e_1+e_3)$ are noncomparable $A$-invariant faces of
$\IR^{\,\,3}_{\,\,+}$ associated with $1$.
So,
in this case,
condition (a) of Theorem 7.2 is satisfied,
but not condition (f).
Example 5.5 can be used to show that the implication (f)$\Longrightarrow$(a)
also does not hold.
It is not clear to us whether condition (f) is sufficient for condition (b) of
Theorem 7.2. \par\parindent=16 pt

Now we give some further remarks in concern with conditions (a)--(e) of Theorem
7.2.
\par\parindent=16 pt

When $\lambda=\rho(A)$ and $A$ is a nonnegative matrix,
by Theorem 7.2 (iii) \par\parindent=0 pt conditions (a)--(d) are all
equivalent.
However,
condition (e) is a strictly weaker condition.
To see this,
consider the following $4\!\times\! 4$ nilpotent \par nonnegative matrix:
\begin{center}
$A=
\left[
\begin{array}{cccc}
0 & 1 & 1 & 0 \\
0 & 0 & 0 & 1 \\
0 & 0 & 0 & 1 \\
0 & 0 & 0 & 0
\end{array}
\right]
$
.
\end{center}
\par\parindent=0 pt
In this case,
each of the singletons $\{ 1\}, \{ 2\}, \{ 3\}$ and $\{ 4\}$ is a class,
indeed a \par\parindent=0 pt basic class,
of $A$.
Also,
$\{ 1\}$ is the only distinguished class of $A$ (associated with $0$).
By the Frobenius-Victory theorem,
we have $\dim ({\cal N}(A)\bigcap\IR^{\,\,n}_{\,\,+})=1$;
that is,
condition (e) is satisfied for $\lambda=0$.
(In fact,
in this case,
for a similar reason,
the dual condition that $\dim({\cal N}(A^T)\bigcap\IR^{\,\,n}_{\,\,+})=1$ is
also satisfied.)\par
However,
condition (a) is not satisfied for $\lambda=0$,
because $\dim {\cal N}(A)=$\par$4-\rank(A)=2$.
\par\parindent=16 pt

In Section 8 we shall give examples which show that the missing\par\parindent=0
pt
implications of Theorem 7.2 all do not hold in general.
\par\parindent=16 pt

The following example shows that when $\lambda<\rho(A)$,
condition (b) of \par\parindent=0 pt Theorem 7.2 and the corresponding
condition with $A$
replaced by $A^T$ and $K$ replaced by $K^*$ together do not
imply condition (a),
not even for the nonnegative matrix case,
assuming that $\lambda$ is a distinguished eigenvalue of $A$ for $K$ and also
of $A^T$ for $K^*$.
\par\parindent=0 pt
\vskip 12 pt

{\bf Example 7.3.}$\,\,$
Consider the following $4\!\times\! 4$ nonnegative matrix:
\begin{center}
$A=
\left[
\begin{array}{cccc}
0 & 0 & 0 & 0 \\
1 & 1 & 1 & 0 \\
1 & 1 & 2 & 0 \\
1 & 1 & 1 & 0
\end{array}
\right]
$
.
\end{center}
\par\parindent=0 pt
In this case $A$ has precisely three classes which form a chain,
namely,
$\{ 4\}, \{ 2, 3\}$ and $\{ 1\}$;
$\{ 2, 3\}$ being a basic class,
and $\{ 1\}$ and $\{ 4\}$ being nonbasic classes both associated with $0$.
Clearly $\IR^{\,\,4}_{\,\,+}$ has precisely three nonzero $A$-invariant faces,
namely,
$\Phi(e_4)$,
$\Phi(e_2+e_3+e_4)$ and $\IR^{\,\,4}_{\,\,+}$ itself.
Also,
$\Phi(e_4)$ is the only $A$-invariant face of
$\IR^{\,\,4}_{\,\,+}$ associated with $0$,
indeed,
a semi-distinguished $A$-invariant face;
thus condition (b) of Theorem 7.2 is satisfied by $A$ for $\lambda=0$.
For a similar reason,
$A^T$ also satisfies the corresponding condition.
However,
condition (a) of Theorem 7.2 is not satisfied,
because we have
$\dim {\cal N}(A)=4-\rank A=2$.
\par\parindent=16 pt
\vskip 12 pt

In the nonpolyhedral case,
even when $K$ is the $n$-dimensional ice-cream cone $K_n$
(which is a perfect cone,
and possesses many nice properties as described in Theorem 6.6)
and $\lambda=\rho(A)$,
condition (b) of Theorem 7.2 and the corresponding condition for $A^T$ together
do not guarantee
condition (a).
For instance,
take $A$ to be $yz^T$,
where $y$ and $z$ are any two mutually orthogonal nonzero vectors that
lie in $\partial K_n$.
Then as can be readily seen,
$\Phi(y)$ and $K_n$ itself are the only nonzero $A$-invariant faces of $K_n$,
both being semi-distinguished $A$-invariant.
So condition (b) of Theorem 7.2 is satisfied.
For a similar reason,
$A^T$ also satisfies the corresponding condition.
However,
$A$ does not satisfy condition (a),
because
\begin{center}
$\dim{\cal N}(\rho(A)I-A)=\dim{\cal N}(A)=n-\rank A=n-1\ge 2$.
\end{center}
\par\parindent=16 pt

Conceivably,
the geometry of $K$ affects the spectral properties of its
cone-preserving maps.
We just mention two such results on a rank-one cone-preserving map,
whose proofs are not difficult.
\par\parindent=0 pt
\vskip 12 pt

{\bf Theorem 7.4.}
{\it Given ${\bf 0}\not= y\in K$ and ${\bf 0}\not= z\in K^*$ such that
$z^Ty=0$.
Let $A=yz^T$.
Then for $\lambda=\rho(A)$,
conditions} (b),
(c) {\it and} (d) {\it of Theorem 7.2 are equivalent,
and the following is another equivalent condition:}
\par\parindent=16 pt
(a)$'$ {\it $y$ is an exposed extreme vector of $K$ such that
$\Phi(y)=d_{K^*}(\Phi(z))$
and the set of all faces of $K$ that contain $y$ forms a chain.}
\par\parindent=0 pt
\vskip 12 pt

{\bf Theorem 7.5.}
{\it Let $y$ and $z$ be nonzero vectors of $\partial K$ and $\partial K^*$
respectively,
and let $A=yz^T$.
The following conditions are equivalent:}
\par\parindent=16 pt
(a) {\it $K$ is an $A$-invariant join-irreducible face of itself.}
\par\parindent=16 pt
(b) {\it $y$ belongs to precisely one maximal face of $K$,
which also includes $d_{K^*}(\Phi(z))$.}
\par\parindent=0 pt
{\it If,
in addition,
the vectors $y$ and $z$ are orthogonal,
then the following is another equivalent condition:}
\par\parindent=16 pt
(c) {\it $y$ belongs to precisely one maximal face of $K$.}
\par\parindent=0 pt
\vskip 25 pt

\begin{center}
8. F{\scriptsize URTHER EXAMPLES}
\end{center}
\par\parindent=0 pt

{\bf Example 8.1.}$\,\,$
\par\parindent=0 pt
\vskip 6 pt
Let $A=
\left[
\begin{array}{cccc}
1 & 1 & 0 & 0	   \\
0 & 1 & 0 & 0	   \\
0 & 0 & 1 & 0	   \\
0 & 0 & 0 & \alpha
\end{array}
\right]
$
,
where $0<\alpha <1$.
\par\parindent=0 pt
\vskip 6 pt
Let $K_1, K_2$ be convex cones of $\IR^{\,\,4}$ given by:
\hskip 50 pt
\begin{tabbing}
\qquad\qquad\qquad
\=$K_1=\pos\{ e_1, A^i(e_2-e_4), i=0,1,\cdots\}$,\\
\> and$\,\,\,$ \' $K_2=\pos\{ e_1, A^i(e_2+e_3+e_4), i=0,1,\cdots\}$,
\end{tabbing}
\par\parindent=0 pt
where $e_1,\cdots,e_4$ are the standard unit vectors of $\IR^{\,\,4}$.
Also let $K=K_1+K_2$.
Clearly $K_1$,
$K_2$ are closed,
pointed cones,
and $K$ is a full cone of $\IR^{\,\,4}$.
As can be easily seen,
$K_1\bigcap (-K_2)=\{ {\bf 0}\}$.
So by [Roc, Corollary 9.1.2] $K$ is a closed,
pointed,
full,
and hence a proper cone of $\IR^{\,\,4}$.
Since $K_1$ and $K_2$ are both invariant under $A$,
we have $A\in\pi(K)$.
Note that the vector $e_1$ (and its nonnegative multiples) is the only common
vector of $K_1$ and $K_2$.
Except for the vector $e_1$,
the second component of each vector in $K_1$ is always positive and the third
component is always zero.
On the other hand,
the second and the third component of each vector in $K_2$,
except for the vector $e_1$,
are always equal and positive.
With these observations in mind,
it is straightforward to verify that $K_1$ and $K_2$ are both faces of $K$.
\par\parindent=16 pt

Clearly each of the faces $\Phi(e_1), K_1, K_2$ and $K$ itself are
$A$-invariant.
\par\parindent=0 pt Indeed,
they are precisely all the nonzero $A$-invariant faces of $K$.
To see this,
let $F$ be any nonzero $A$-invariant face of $K$.
If $F$ contains an extreme vector of the form $A^i(e_2-e_4)$ for some $i\ge 0$,
then $F$ will contain all the subsequenct vectors $A^k(e_2-e_4)$ for
$k=i+1,i+2,\cdots$
and hence a relative interior vector of $K_1$;
thus $F$ includes $K_1$.
Similarly,
$F$ includes $K_2$,
if $F$ contains an extreme vector of the form $A^i(e_2+e_3+e_4)$ for some $i\ge
0$.
But the extreme vectors of $K$
are (up to multiples) precisely the vectors $e_1$ and $A^i(e_2-e_4)$,
$A^i(e_2+e_3+e_4)$,
$i=0,1,2,\cdots$.
With these observations,
by considering the possible extreme vectors that $F$ may have,
we readily see that $F$ must be one of the faces $\Phi(e_1), K_1, K_2$ or $K$.
\par\parindent=16 pt

Note that $K$ is $A$-invariant join-reducible,
as $K=K_1\bigvee K_2$.
Also each of the faces $K_1, K_2$ does not contain in its relative interior
a generalized eigenvector of $A$.
So,
$\Phi(e_1)$ is the only semi-distinguished $A$-invariant face of $K$,
and hence,
with $\lambda=\rho(A)$,
condition (b) of Theorem 7.1 is satisfied trivially.
However,
as we are going to show,
in this case we have $m_{\rho(A)}=2$,
so that condition (a) of Theorem 7.1 is not satisfied.
Thus,
this example shows that in general
conditions (a) and (b) of Theorem 7.1 are not equivalent.
\par\parindent=16 pt

Clearly,
the vector $e_2-e_4+A^2(e_2-e_4)$ ($={(2,2,0,-(1+\alpha^2))}^T$) belongs
\par\parindent=0 pt to $\relint K_1$.
Similarly,
the vector $e_2+e_3+e_4+A^2(e_2+e_3+e_4)$ ($={(2,2,2,1+\alpha^2)}^T$) belongs
to $\relint K_2$.
Hence,
their sum,
which is ${(4,4,2,0)}^T$,
belongs to $\int K$.
The latter vector is,
in fact,
also a generalized eigenvector of $A$ corresponding to $\rho(A)$
of order two.
This shows that
$m_{\rho(A)}=2$,
as desired.
\par\parindent=16 pt

Since $\Phi(e_1)$ and $K$ are the only $A$-invariant faces that contain
in their relative interiors a generalized eigenvector of $A$,
with $\lambda=\rho(A)$,
condition (c) of Theorem 7.2 is satisfied.
However,
condition (b) of Theorem 7.2 is not satisfied (with $\lambda=\rho(A)$),
as $K$ is not semi-distinguished $A$-invariant.
Thus,
this example shows that,
in general,
for the conditions (b),
(c) of Theorem 7.2,
the implication (c)$\Longrightarrow$(b) does not hold.
\par\parindent=16 pt

Also note that in this example there does not exist a semi-distinguished
$A$-invariant face $F$ such that $\sp_A(F)=\sp_A(K)$,
though $K$ contains in its interior a generalized eigenvector of $A$.
So the converse of Remark 6.5 is not true.
\par\parindent=0 pt
\vskip 12 pt

{\bf Example 8.2.}$\,\,$
Take any fixed positive real number $\alpha$ less than $1$,
and let
\begin{center}
$A=
\left[
\begin{array}{cccc}
1 & 1 & 0 & 0	   \\
0 & 1 & 0 & 0	   \\
0 & 0 & 1 & 0	   \\
0 & 0 & 0 & \alpha
\end{array}
\right]
$
.
\end{center}
\par\parindent=0 pt
Let $K_1, K_2$ be the convex cones of $\IR^{\,\,4}$ given by:
\begin{tabbing}
\qquad\qquad\qquad
\=$K_1=\pos\{ e_1, e_4, A^i(e_2-e_4), i=0,1,2,\cdots\}$,\\
\> and$\,\,\,$ \' $K_2=\pos\{ e_1, e_4, A^i(e_2+e_3+e_4), i=0,1,2,\cdots\}$.
\end{tabbing}
\par\parindent=0 pt
Also,
let $K=K_1+K_2$.
Then $K$ is a proper cone of $\IR^{\,\,4}$ and we have $A\in\pi(K)$.
(Note that the cones $K_1, K_2$ considered here are each obtained from the
corresponding
cones in Example 8.1 by adjoining the extra extreme ray $\Phi(e_4)$,
whereas the matrix $A$ is the same as before.)
\par\parindent=16 pt

Using an argument similar to the one used in Example 8.1,
one can readily show that the nonzero $A$-invariant faces of $K$ are precisely
$\Phi(e_1), \Phi(e_4),$\par\parindent=0 pt $\Phi(e_1+e_4), K_1, K_2$ and $K$
itself.
Among these faces,
only the faces $\Phi(e_1), K_1$ and $K$ have the property that it contains in
its
relative interior a \par generalized eigenvector of $A$ corresponding to
$\rho(A)$.
Furthermore,
$\Phi(e_1)$ and $K_1$ are the only semi-distinguished $A$-invariant faces
associated with $\rho(A)$,
and $\sp_A(K_1)=\sp_A(K)=(1, 2)$.
Thus,
with $\lambda=\rho(A)$,
condition (c) of \par Theorem 6.4 is satisfied.
However,
condition (a)(i) of Theorem 6.4 is not \par satisfied,
as $K$ cannot be written as the join of semi-distinguished $A$-invariant
faces associated with $\rho(A)$.
Also note that $K_2$ is a nonzero $A$-invariant face \par\parindent=0 pt
associated with $\rho(A)$,
but there does not exist a semi-distinguished $A$-invariant face $F$ included
in $K_2$ such that $\sp_A(F)=\sp_A(K_2)=(1, 2)$.
So condition (b)(i) of Theorem 6.4 is also not fulfilled.
Thus,
this example shows that,
in general,
condition (c) of Theorem 6.4 does not imply condition (a) nor condition (b)
there.
\par\parindent=0 pt
\vskip 12 pt

{\bf Example 8.3.}$\,\,$
Let $A=J_2(1)\bigoplus J_1(1)\bigoplus J_1(\alpha)\bigoplus J_1(1)$,
where $0<\alpha <1$.\par\parindent=0 pt
On $\IR^{\,\,5}$ define the following cones:
\begin{tabbing}
\qquad\qquad\qquad
\= $K_1=\pos\{ e_1, A^i(e_2-e_4), i=0,1,2,\cdots\}$,\\
\> $K_2=\pos\{ e_1, A^i(e_2+e_3+e_4), i=0,1,2,\cdots\}$,\\
\> $K_3=\pos\{ e_1, A^i(e_2+e_4+e_5), i=0,1,2,\cdots\}$,
\end{tabbing}
\par\parindent=0 pt
and $K=K_1+K_2+K_3$.
Clearly $K_1, K_2$ and $K_3$ are each closed,
pointed cones.
Furthermore,
it is easy to verify that
\begin{center}
$K_1\bigcap (-K_2)=\{ {\bf 0}\}\,\,$ and $\,\, (K_1+K_2)\bigcap (-K_3)=\{ {\bf
0}\}$.
\end{center}
\par\parindent=0 pt
It follows that $K$ is a proper cone of $\IR^{\,\,5}$.
\par\parindent=16 pt

It is not difficult to show that $K_1, K_2$ and $K_3$ are all faces of $K$.
In addition,
we can also show that $K_1+K_2$ and $K_1+K_3$ are also faces of $K$.
(Whether $K_2+K_3$ is a face of $K$ is not obvious to us,
but we do not need that for our purposes.)
Then clearly $K_1+K_2=K_1\bigvee K_2$ and $K_1+K_3=K_1\bigvee K_3$.
\par\parindent=16 pt

Clearly,
we have $A\in\pi(K)$.
It is not difficult to show that
the nonzero $A$-invariant faces of $K$ are precisely
$\Phi(e_1), K_1, K_2, K_3,K_1+K_2, K_1+K_3, K$,
and possibly $K_2+K_3$.
Among these $A$-invariant faces of $K$,
only $\Phi(e_1)$ is semi-distinguished $A$-invariant.
Note that $K_1+K_2, K_1+K_3$ are \par\parindent=0 pt noncomparable
$A$-invariant faces,
each containing in its relative interior a generalized eigenvector of $A$
corresponding to $\rho(A)$.
So,
in this case,
with $\lambda=\rho(A)$,
condition (d) of Theorem 7.2 is satisfied,
but not condition (c).
This shows that for conditions (c),
(d) of Theorem 7.2,
in general,
we do not have (d)$\Longrightarrow$(c).
\par\parindent=0 pt
\vskip 25 pt
\newpage

\begin{center}
9. O{\scriptsize PEN QUESTIONS}
\end{center}
\par\parindent=16 pt

In trying to improve or better understand the result of Theorem 6.6,
we pose the following question.
\par\parindent=0 pt
\vskip 12 pt

{\it Question} 9.1.$\,\,$
Let $K$ be a proper cone whose dual cone $K^*$ is a facially exposed cone.
Is it true that for any $A\in\pi(K)$,
we have the following?
\par\parindent=16 pt
(i)
For any nonzero $A$-invariant face $F$ of $K$,
$F$ is semi-distinguished $A$-invariant if and only if $\sp_A(F)\succ\sp_A(G)$
for all $A$-invariant faces $G$ properly included in $F$.
\par\parindent=16 pt
(ii)
There exists in $K$ a generalized eigenvector of $A$ corresponding to $\rho(A)$
of order $\nu_{\rho(A)}(A)$.
\par\parindent=16 pt
(iii)
For any $A$-invariant face $F$ which contains in its relative interior a
generalized eigenvector of $A$,
there exists a semi-distinguished $A$-invariant face $G$ included in $F$ such
that
$\sp_A(G)=\sp_A(F)$.
\par\parindent=0 pt
\vskip 12 pt

{\it Question} 9.2.$\,\,$
Let $K$ be a proper cone with the property that the dual cone of each of its
faces is a facially exposed cone.
Is it true that,
for any $A\in\pi(K)$ and any distinguished eigenvalue $\lambda$ of $A$,
condition (a) of Theorem 6.4 is always satisfied?
\par\parindent=16 pt
\vskip 12 pt

We do not know the answer to Question 9.2 even when $K$ is a perfect cone or is
equal to $P(n)$ for some positive integer $n$.
If the answer is in the negative,
then it will follow that conditions (a) and (b) of Theorem 6.4 are logically
independent
(cf.
Theorem 6.6 and the discussion following Remark 6.5).
\par\parindent=0 pt
\vskip 25 pt
\newpage

\begin{center}
R{\scriptsize EFERENCES}
\end{center}
\par\parindent=16 pt

\footnotesize
{
\par\parindent=0 pt
$[$Bar 1$]$ $\,\,$
G.P. Barker,
{\it Perfect cones},
Linear Algebra Appl. {\bf 22} (1978),
211--221.
\vskip 5 pt
\par\parindent=0 pt
$[$Bar 2$]$ $\,\,$
$\lline$,
{\it Theory of cones},
Linear Algebra Appl. {\bf 39} (1981),
263--291.
\vskip 5 pt
\par\parindent=0 pt
$[$B--P$]$ $\,\,$
A. Berman
and R.J. Plemmons,
{\it Nonnegative Matrices in the Mathematical \par\hskip 30 pt Sciences},
Academic Press,
New York,
1979.
\vskip 5 pt
\par\parindent=0 pt
$[$B--T$]$ $\,\,$
G.P. Barker
and A. Thompson,
{\it Cones of polynomials},
Portugal. Math. {\bf 44} (1987), \par\hskip 30 pt
183--197.
\vskip 5 pt
\par\parindent=0 pt
$[$Bir$]$ $\,\,$
G. Birkhoff,
{\it Lattice Theory},
3rd ed.,
American Mathematical Society,
Providence, \par\hskip 30 pt
R.I.,
1966.
\vskip 5 pt
\par\parindent=0 pt
$[$Dod$]$ $\,\,$
P.G. Dodds,
{\it Positive compact operators},
Quaestiones Math. {\bf 18} (1995),
21--45.
\vskip 5 pt
\par\parindent=0 pt
$[$G--L--R$]$ $\,\,$
L. Gohberg,
P. Lancaster,
and L. Rodman,
{\it Invariant Subspaces of Matrices \par\hskip 30 pt with Applications},
John Wiley \& Sons,
New York,
1986.
\vskip 5 pt
\par\parindent=0 pt
$[$H--S$]$ $\,\,$
D. Hershkowitz
and H. Schneider,
{\it On the generalized nullspace of $M$-matrices and \par\hskip 30 pt
$Z$-matrices},
Linear Algebra Appl. {\bf 106} (1988),
5--23.
\vskip 5 pt
\par\parindent=0 pt
$[$J--V 1$]$ $\,\,$
R. Jang
and H.D. Victory, Jr.,
{\it On nonnegative solvability of linear integral \par\hskip 30 pt equations},
Linear Algebra Appl. {\bf 165} (1992),
197--228.
\vskip 5 pt
\par\parindent=0 pt
$[$J--V 2$]$ $\,\,$
$\lline$,
{\it On the ideal structure of positive,
eventually compact linear operators on \par\hskip 30 pt Banach lattices},
Pacific J. Math. {\bf 157} (1993),
57--85.
\vskip 5 pt
\par\parindent=0 pt
$[$J--V 3$]$ $\,\,$
$\lline$,
{\it On nonnegative solvability of linear operator equations},
Integr. Equat. \par\hskip 30 pt Oper. Th. {\bf 18} (1994),
88--108.
\vskip 5 pt
\par\parindent=0 pt
$[$Nel$]$ $\,\,$
P. Nelson, Jr.,
{\it The structure of a positive linear integral operator},
J. London Math. \par\hskip 30 pt Soc. (2) {\bf 8} (1974),
711--718.
\vskip 5 pt
\par\parindent=0 pt
$[$MN 1$]$ $\,\,$
P. Meyer-Nieberg,
{\it A partial spectral reduction for positive linear operators},
Arch. \par\hskip 30 pt Math. {\bf 45} (1985),
34--41.
\vskip 5 pt
\par\parindent=0 pt
$[$MN 2$]$ $\,\,$
$\lline$,
{\it Banach Lattices},
Springer-Verlag,
New York,
1991.
\vskip 5 pt
\par\parindent=0 pt
$[$Roc$]$ $\,\,$
R.T. Rockafellar,
{\it Convex Analysis},
Princeton Univ. Press,
Princeton,
NJ,
1970.
\vskip 5 pt
\par\parindent=0 pt
$[$Rot$]$ $\,\,$
U.G. Rothblum,
{\it Algebraic eigenspaces of nonnegative matrices},
Linear Algebra \par\hskip 30 pt Appl. {\bf 12} (1975),
281--292.
\vskip 5 pt
\par\parindent=0 pt
$[$Su--T$]$ $\,\,$
C.H. Sung
and B.S. Tam,
{\it A study of projectionally exposed cones},
Linear Algebra \par\hskip 30 pt Appl. {\bf 139} (1990),
225--252.
\vskip 5 pt
\par\parindent=0 pt
$[$Scha$]$ $\,\,$
H.H. Schaefer,
{\it Banach lattices and positive operators},
Springer Verlag,
Berlin-\par\hskip 30 pt Heidelberg-New York,
1974.
\vskip 5 pt
\par\parindent=0 pt
$[$Schn 1$]$ $\,\,$
H. Schneider,
{\it The elementary divisors associated with $0$ of a singular $M$-matrix},
\par\hskip 30 pt
Proc. Edinburgh Math. Soc. (2) {\bf 10} (1956),
108--122.
\vskip 5 pt
\par\parindent=0 pt
$[$Schn 2$]$ $\,\,$
$\lline$,
{\it Geometric conditions for the existence of positive eigenvalues of
matrices},
\par\hskip 30 pt
Linear Algebra Appl. {\bf 38} (1981),
253--271.
\vskip 5 pt
\par\parindent=0 pt
$[$Schn 3$]$ $\,\,$
$\lline$,
{\it The influence of the marked reduced graph of a nonnegative matrix on
\par\hskip 30 pt the Jordan
form and on related properties$\,$}:
{\it a survey},
Linear Algebra Appl. {\bf 84} \par\hskip 30 pt (1986),
161--189.
\vskip 5 pt
\par\parindent=0 pt
$[$Tam 1$]$ $\,\,$
B.S. Tam,
{\it A note on polyhedral cones},
J. Austral. Math. Soc. Ser A {\bf 22} (1976), \par\hskip 30 pt
456--461.
\vskip 5 pt
\par\parindent=0 pt
$[$Tam 2$]$ $\,\,$
$\lline$,
{\it On the duality operator of a convex cone},
Linear Algebra Appl. {\bf 64} \par\hskip 30 pt (1985),
33--56.
\vskip 5 pt
\par\parindent=0 pt
$[$Tam 3$]$ $\,\,$
$\lline$,
{\it On the distinguished eigenvalues of a cone-preserving map},
Linear \par\hskip 30 pt Algebra Appl. {\bf 131} (1990),
17--37.
\vskip 5 pt
\par\parindent=0 pt
$[$Tam 4$]$ $\,\,$
$\lline$,
{\it On semipositive bases for a cone-preserving map},
in preparation.
\vskip 5 pt
\par\parindent=0 pt
$[$T--S 1$]$ $\,\,$
B.S. Tam
and H. Schneider,
{\it On the core of a cone-preserving map},
Trans. Am. \par\hskip 30 pt Math. Soc. {\bf 343} (1994),
479--524.
\vskip 5 pt
\par\parindent=0 pt
$[$T--S 2$]$ $\,\,$
$\lline$,
{\it Solutions of linear equations over cones},
in preparation.
\vskip 5 pt
\par\parindent=0 pt
$[$T--W$]$ $\,\,$
B.S. Tam
and S.F. Wu,
{\it On the Collatz-Wielandt sets associated with a cone-\par\hskip 30 pt
preserving map},
Linear Algebra Appl. {\bf 125} (1989),
77--95.
\vskip 5 pt
\par\parindent=0 pt
$[$Vic 1$]$ $\,\,$
H.D. Victory, Jr.,
{\it On linear integral operators with nonnegative kernels},
J. Math. \par\hskip 30 pt Anal. Appl. {\bf 89} (1982),
420--441.
\vskip 5 pt
\par\parindent=0 pt
$[$Vic 2$]$ $\,\,$
$\lline$,
{\it The structure of the algebraic eigenspace to the spectral radius of
eventually
\par\hskip 30 pt
compact,
nonnegative integral operators},
J. Math. Anal. Appl. {\bf 90} (1982),
484--516.
\vskip 5 pt
\par\parindent=0 pt
$[$Zer$]$ $\,\,$
M. Zerner,
{\it Quelques propri${\acute {e}}$t${\acute {e}}$s spectrales des op${\acute
{e}}$rateurs positifs},
J. Funct. Anal. \par\hskip 30 pt {\bf 72} (1987),
381--417.

}

\par\parindent=20 pt
\vskip 20 pt

{\small D}{\scriptsize EPARTMENT OF} {\small M}{\scriptsize ATHEMATICS},
{\small T}{\scriptsize AMKANG} {\small U}{\scriptsize NIVERSITY},
{\small T}{\scriptsize AMSUI},
{\small T}{\scriptsize AIWAN} {\small 25137},
{\small ROC}
\par\parindent=20 pt
{\small\it E-mail address} :
{\small bsm01{$@$}mail.tku.edu.tw}
\vskip 10 pt
\par\parindent=20 pt

{\small D}{\scriptsize EPARTMENT OF} {\small M}{\scriptsize ATHEMATICS},
{\small U}{\scriptsize NIVERSITY OF} {\small W}{\scriptsize ISCONSIN},
{\small M}{\scriptsize ADISON},
{\small W}{\scriptsize ISCONSIN} {\small 53706},
{\small USA}
\par\parindent=20 pt
{\small\it E-mail address} :
{\small hans{$@$}math.wisc.edu}

}

\end{document}